\documentclass[11pt]{article}
\usepackage{amscd}
\usepackage{amssymb}
\newtheorem{lemma}{Lemma}
\newtheorem{theorem}{Theorem}
\newtheorem{corollary}{Corollary}
\newtheorem{definition}{Definition}
\newtheorem{prop}{Proposition}
\newtheorem{step}{Step}

\newcommand{\C}{\mathbb{C}}
\newcommand{\R}{\mathbb{R}}

\newcommand{\Proj}{\mathbb{P}}
\newcommand{\qed}{\nopagebreak\begin{flushright}\framebox[2mm]{}\end{flushright}\vspace{5 mm}}
\title{Deformations of Asymptotically Conical Special Lagrangian Submanifolds}
\date{}
\author{Tommaso Pacini}
\begin{document}
\maketitle
\textbf{Abstract:} McLean [ML] studied the deformations of compact special Lagrangian submanifolds, showing in particular that they come in smooth moduli spaces whose dimension depends only on the topology of the submanifold. In this article we study the analogous problem for non-compact, ``asymptotically conical'' SL submanifolds, with respect to various ``boundary conditions at infinity''.
\section{Introduction}
In the context of Riemannian geometry, much research is directed to the study of minimal submanifolds. One of the questions that arises naturally in this field is the following.

\ 

\textbf{Question:} Can a given minimal submanifold be ``deformed'' to get new examples? What parameters might be involved in these deformations?

\ 

The naive approach is to parametrize these deformations as the zero-set of a ``mean curvature operator'', then study them using the implicit function theorem. However, this entails a good understanding of the Jacobi operator of the initial submanifold $\Sigma$, which in general is not possible.

The work of Oh and, more recently, of McLean (cfr. [Oh], [ML]) shows that, in the ``right'' geometric context, the problem simplifies and sometimes becomes tractable. We are thus lead to the study of minimal Lagrangian submanifolds in Kaehler-Einstein (KE) ambient spaces, and of special Lagrangian (SL) submanifolds in Calabi-Yau (CY) manifolds: cfr. section \ref{section_SL} for definitions and examples.

In [ML], McLean shows that every ``infinitesimal SL deformation'' of a (smooth) compact SL submanifold is ``integrable''; i.e., it generates actual SL deformations. A corollary of this is, in the compact case, that the set
$$\mathcal{M}(\Sigma):=\{\mbox{SL immersions }\Sigma\hookrightarrow N\}$$
has a natural smooth structure. The dimension of this ``moduli space'' turns out to be related only to the topology of $\Sigma$.

For at least two reasons, it is interesting to understand if similar results hold also for non-compact SL submanifolds.

The first (perhaps simplistic) reason is that it is well-known that our basic ``model'' CY manifold, $\C^n$ with its standard structures, does not admit compact minimal submanifolds. In particular, it cannot contain compact SL submanifolds.

The second reason is more interesting. Moduli spaces of smooth, compact, SL submanifolds are not, in general, compact; it is however expected that they may be compactified by adding ``boundaries'' consisting of singular SLs. 

Linear algebra, based at any point $p$ of a CY manifold $N$, shows that $T_pN$, with all its structures, is isomorphic to $\C^n$. The simplest type of singularity of a SL submanifold is given by isolated points $p$ whose ``tangent space'' in $T_pN$ is isomorphic to a SL cone in $\C^n$. The submanifold itself should locally be the limit of smooth ``asymptotically conical'' (AC) SL submanifolds. Understanding the class of AC SL submanifolds should thus be a fundamental element both in compactification and in desingularization procedures (such as gluing).

McLean's result relies on the standard Hodge theory for compact manifolds. In the AC case, the basic tool is provided by the Fredholm theory of the Laplace operator, due to Nirenberg, Walker, Cantor, McOwen, Lockhart, Melrose, Bartnik et al.; cfr. in particular [NW], [MO], [B], [M1], [M2]. 

The most complete set of results in this direction has been developed by R. Melrose. In section \ref{section_ACintro} we thus introduce asymptotically conical and asymptotically cylindrical manifolds using the language of sc- and b-geometry, following [M1] and [M2]. 

In section \ref{section_laplace}, still following [M1], we present various results regarding intrinsic analysis on these manifolds (referring elsewhere for detailed proofs whenever the statements or techniques seem sufficiently well-known). In section \ref{subsection_1forms}, we then focus on certain harmonic 1-forms, corresponding to the ``infinitesimal AC SL deformations'': we prove a representation theorem and calculate the dimensions of these spaces, with respect to the two most interesting ``rates of decay at infinity''.

In section \ref{section_immersions} we discuss the ``AC'' condition for submanifolds, showing how to preserve it under normal deformations. This requires a study of the exponential map and Jacobi field estimates, which may be of independent interest in b-geometry.

Section \ref{section_SLdefs} presents the main result: as in the compact case, every infinitesimal AC SL deformation is integrable; i.e., each such inf. deformation generates a 1-parameter family of SL submanifolds which, by the results of section \ref{section_immersions}, are also AC. The results from section \ref{subsection_1forms} then show that the dimension of the space of deformations depends both on the topology and and on the analytic properties of $(\Sigma,g)$. 

We also discuss the role played by the curvature of the ambient space in the construction of spaces of deformations having different speeds of decay (eg: $L^2$ decay) at infinity.

\textbf{Note:} Throughout the article, we will assume $n>2$.

\textbf{Acknowledgements:} I wish to thank P. de Bartolomeis and G. Tian for their support, suggestions and interest. I am also grateful to R. Melrose for many useful explanations on scattering geometry and to A. Ghigi and P. Piazza for interesting conversations.

The problem of understanding the deformations of asymptotically conical SL submanifolds was posed by D. Joyce in a previous version of [J1]. This work started as an attempt to understand that article. A similar result to that presented here was more recently conjectured in the updated version of [J1] and in [J2]. After the present work was completed, but independently from it, Joyce's conjectures were proved by S. Marshall (cfr. [M]). For further comments, see remark 5, section \ref{section_SLdefs}.

Finally, I gratefully acknowledge the generous support of the University of Pisa and of GNSAGA, and the hospitality of MIT.
\section{Minimal Lagrangian and special Lagrangian submanifolds} \label{section_SL}
Recall the following, standard, definitions.
\begin{definition} An immersion $\phi:\Sigma\hookrightarrow N$ of a manifold $\Sigma$ into a Riemannian manifold $(N,g)$ is ``minimal'' if the corresponding mean curvature vector field vanishes: $H_\phi\equiv 0$.

An immersion $\phi:\Sigma\hookrightarrow N$ into a symplectic manifold $(N,\omega)$ is ``isotropic'' if $\phi^*\omega\equiv 0$; if dim $\Sigma=n$ and dim $N=2n$, isotropic submanifolds are called ``Lagrangian''.
\end{definition}
Usually the normal bundle $(T\Sigma)^\perp$ of a submanifold $\Sigma\subseteq (N,g)$ is not an intrinsic object: it depends on the immersion and on the ambient space.

However, assume that the ambient manifold is a Kaehler manifold: $(N^{2n},J,g,\omega)$, where $\omega(X,Y):=g(JX,Y)$ is the Kaehler 2-form. We will use the following notation:
\begin{itemize}
\item $K_N:=\Lambda^{n,0}(N)$ is the canonical bundle of $(N,J)$;
\item $\nabla$ and $Ric$ are the Levi-Civita connection and the symmetric Ricci 2-tensor associated to $(N,g)$;
\item $\rho(X,Y):=Ric(JX,Y)$ is the ``Ricci 2-form''.
\end{itemize}
\begin{definition} A Kaehler manifold $(N,J,g,\omega)$ is ``Kaehler-Einstein'' if $Ric=c\cdot g$, for some $c\in \R$; equivalently, if $\rho=c\cdot \omega$.
\end{definition}
If $\Sigma$ is Lagrangian and $N$ is Kaehler, the (restriction of the) symplectic form gives an isomorphism 
$$(T\Sigma)^\perp\simeq\Lambda^1(\Sigma),\ \ \ \ V\simeq\nu:=\omega(V,\cdot)_{|\Sigma}$$
This shows that $(T\Sigma)^\perp$ is actually independent of the ambient manifold $N$.

We now want to introduce ``Calabi-Yau'' manifolds and ``special Lagrangian'' submanifolds.

Given a smooth section $\Omega$ of $K_N$, recall that $d\,\Omega=(\partial+\overline{\partial})\Omega=\overline{\partial}\,\Omega$, since $\partial \,\Omega\in \Lambda^{n+1,0}(N)\equiv 0$. Thus $\Omega$ is closed iff $\Omega $ is holomorphic.

Furthermore, $\nabla\Omega\equiv 0\Rightarrow d\,\Omega\equiv 0$.

\begin{definition} A (differentiable) manifold $N^{2n}$ is of ``Calabi-Yau (CY) type'' if it admits a Riemannian metric $g$ with holonomy $Hol(N,g)\leq SU(n;\C)$; equivalently, if it admits a Kaehler structure $(J,g,\omega)$ and a (non-zero) section $\Omega$ of $K_N$ such that $\nabla\Omega\equiv 0$.

A choice of $(J,g,\omega,\Omega)$ defines a ``Calabi-Yau structure''.
\end{definition}
In particular, CY manifolds are KE manifolds with $c=0$, i.e. $Ric\equiv 0$.

Since $\Omega$ is parallel, it depends only on some $\Omega[p]$ and can thus be normalized in such a way that $\|\Omega\|=1$. We will always assume this. Such an $\Omega$ is then unique up to a multiplicative factor $\theta\in S^1$.

Given a CY manifold $(N^{2n},J,g,\Omega)$, we will let $\alpha,\beta$ denote the real and imaginary parts of $\Omega$: $\Omega=\alpha+i\beta$. Notice that $d\,\alpha=d\,\beta=0$.

\begin{definition} Let $(N,\Omega)$ be a CY manifold. An oriented, immersed submanifold $\phi:\Sigma^n \hookrightarrow N^{2n}$ is special Lagrangian (SL) iff $\Omega_{|\Sigma}=vol_{\Sigma}$.
\end{definition}

It is simple to prove (cfr. [HL]) that 
$$\Omega_{|\Sigma}=\pm \,vol_{\Sigma} \Leftrightarrow \phi^*\omega=0, \phi^*\beta=0$$
In other words, up to a change of orientation, $\Sigma$ is SL iff $\phi^*\omega=0=\phi^*\beta$. SL submanifolds are thus Lagrangian.

Furthermore, CY manifolds and SL submanifolds are one of the basic examples of ``calibrated geometry'' (cfr. [HL]): in particular, it can easily be shown that special Lagrangians are minimal.

The following, well-known, proposition states that, in CY manifolds, minimal Lagrangian and SL submanifolds are closely related.

\begin{prop}
Let $(N,J,g,\Omega)$ be a CY manifold and $\Sigma\subseteq N$ be an oriented Lagrangian submanifold. Then $\Sigma$ is minimal iff it is SL with respect to the CY structure $\theta\cdot \Omega$, for some $\theta\in S^1$.
\end{prop}

Finding examples of SL submanifolds is not easy; one reason for this is that most examples of CY ambient spaces are provided by an abstract existence theorem due to Yau (cfr. [Y]). However, a few CY manifolds are known explicitly, allowing for explicit constructions of SL submanifolds.

The basic example is $\C^n\simeq \R^{2n}$ with its standard structures $J_{std},g_{std},\omega_{std}$ and $\Omega_{std}:=dz^1\wedge\dots\wedge dz^n$. Many examples of SL submanifolds in $\C^n$ are now known: cfr. eg [HL], [J2], [H], [G]. In particular, [H] shows how minimal Lagrangian submanifolds in $\C\Proj^{n-1}$ give rise to families of SL submanifolds in $\C^n$; an appropriate choice of parametrization then shows that they are ``asymptotically conical'' in the sense of section \ref{section_ACintro}.

A second example of an explicit CY is the following.

Given any compact KE manifold $(N^{2n},J_N,g_N,\omega_N)$ such that $\rho_N=c\cdot \omega_N$, $c>0$, Calabi (cfr. [C]) proved the existence of a complete Kaehler Ricci-flat metric on the canonical bundle $K_N$ of $N$. $K_N$ also has a natural global holomorphic volume form. With respect to these structures, $K_N$ is an ``asymptotically conical'' CY manifold. Again, it is possible to show (cfr. [G], [P]) that minimal Lagrangian submanifolds $\Sigma^n\subseteq N^{2n}$ generate ``asymptotically conical'' SL submanifolds $\hat{\Sigma}^{n+1}\subseteq K_N^{2n+2}$: each $\hat{\Sigma}$ is a $\R$-line bundle over $\Sigma$.

\section{Asymptotically conical and cylindrical manifolds} \label{section_ACintro}
In this section we want to introduce the basic objects of Melrose's ``b-geometry''. Cfr. [M1] for further details.

Let $X^n$ be an oriented compact $n$-dimensional manifold with boundary $\partial X$ and $\Sigma:=X\setminus\partial X$ be its interior. We will use the following notation:

\ 

$2X:=X\bigcup_{\partial X}(-X)$ (compact oriented manifold without boundary).

$C^\infty(X):=$ space of functions on $X$ having a smooth extension on a neighborhood of $X\subseteq 2X$.

$TX:=T(2X)_{|X}$, $T^*X:=T^*(2X)_{|X}$.

If $E$ is a bundle on $X$, $\Lambda^0(E)$ will denote the space of smooth sections and $\Lambda^0_c(E)$ those with compact support in $\Sigma$.

\begin{definition}
$x \in C^\infty(X)$ is a boundary-defining function if $x\geq 0$, $\partial X=\{x=0\}$, $d\,x_{|\partial X}\neq 0$.
\end{definition}

Notice that if $x'$ is a second boundary-defining function, then there exists $a\in C^\infty(X)$, nowhere zero, such that $x'=ax$. Fixing a boundary-defining function $x$ gives a local trivialization of $X$ near $\partial X$: $X\simeq[0,1]\times \partial X$, where the function $x$ corresponds to the coordinate $x$ on $[0,1]$. Now set:

\ 

$\nu_b:=\{Z\in \Lambda^0(TX):Z_{|\partial X} \in T(\partial X)\}$.

$\nu_{sc}:=x\cdot \nu_b=\{Z \in \Lambda^0(TX):Z=xZ',\ Z'\in \nu_b\}$.

\ 

Clearly, $\nu_{sc}$ is independent of the choice of boundary-defining function $x$.

Taylor expansion shows that, locally near $\partial X$, $\nu_b=<x\partial_x,\partial_{y^i}>_{C^\infty(X)}$ and $\nu_{sc}=<x^2\partial_x,x \partial_{y^i}>_{C^\infty(X)}$ where $y^i$ denote local coordinates on $\partial X$.

Thus we have the following chain of spaces of vector fields:
$$\Lambda^0_c(T\Sigma)\subseteq \nu_{sc}\subseteq \nu_b \subseteq \Lambda^0(TX) \subseteq \Lambda^0(T\Sigma)$$
It turns out that one can define vector bundles $^bTX$ and $^{sc}TX$ over $X$, of constant rank $n$, such that $^bTX_{|\Sigma} \simeq T\Sigma \simeq \mbox{$^{sc}TX_{|\Sigma}$}$ and possessing differentiable structures such that $\nu_b=\Lambda^0(^bTX)$,\ $\nu_{sc}=\Lambda^0(^{sc}TX)$.

Let $^bT^*X$, $^{sc}T^*X$ be the corresponding dual bundles. Then, locally near $\partial X$, we may define:

\ 

$\{\frac{d\,x}{x}, d\,y^i\}:=$ the local basis of $^bT^*X$ dual to $\{x\partial_x,\partial_{y^i}\}$.

$\{\frac{d\,x}{x^2},\frac{d\,y^i}{x}\}:=$ the local basis of $^{sc}T^*X$ dual to $\{x^2\partial_x,x\partial_{y^i}\}$. 

\ 

With these definitions, locally $\nu^*_b:=\Lambda^0(^bT^*X)=<\frac{d\,x}{x},d\,y^i>_{C^\infty(X)}$ and $\nu^*_{sc}:=\Lambda^0(^{sc}T^*X)=<\frac{d\,x}{x^2},\frac{d\,y^i}{x}>_{C^\infty(X)}$. Thus:
$$\Lambda^0_c(T^*\Sigma)\subseteq \Lambda^0(T^*X)\subseteq \nu^*_b = x\cdot \nu^*_{sc} \subseteq \nu^*_{sc} \subseteq \Lambda^0(T^*\Sigma)$$
The alternating products of $^bT^*X$, $^{sc}T^*X$ lead to the spaces of b- and sc-differential forms, whose spaces of sections satisfy:
$$\Lambda^i_c(\Sigma) \subseteq \Lambda^i(X)\subseteq \mbox{$^b\Lambda^i(X)$} = x^i\cdot \mbox{$^{sc}\Lambda^i(X)$} \subseteq \mbox{$^{sc}\Lambda^i(X)$} \subseteq \Lambda^i(\Sigma)$$
In this case, we will use the same notation for both the bundles of $i$-forms and their sections. Finally, we set:

\ 

$\mbox{Diff}\,^m_b(X):=$ space of linear operators $P:C^\infty(X) \longrightarrow C^\infty(X)$ such that, locally, 
$$P=\Sigma_{i+|I|\leq m}\,p_{i,I}(x,y)(x\,D_x)^i D_y^I$$ 
where $p_{i,I}\in C^\infty(X)$, $D_x=-\sqrt{-1}\,\partial_x$, $D_{y^i}=-\sqrt{-1}\,\partial_{y^i}$, $I=(i_1,...,i_{n-1})$ is a multi-index, $D_y^I=(D_{y^1})^{i_1}...(D_{y^{n-1}})^{i_{n-1}}$ and $|I|=i_1+...+i_{n-1}$. 

\ 

More generally, for any vector bundles $E$, $F$ over $X$, we define:

\ 

$\mbox{Diff}\,^m_b(E;F):=$ space of linear operators $P:\Lambda^0(E)\longrightarrow \Lambda^0(F)$ such that, locally, all components of $P$ have the above form.

$\mbox{Diff}\,^m_b(E):=\mbox{Diff}\,^m_b(E;E)$.

\

It may be useful to emphasize here that we are working in the category of manifolds with boundary. Thus, when checking that $P\in \mbox{Diff}\,^m_b(E;F)$ near the boundary, one must choose trivializations which hold up to the boundary (clearly, the classes $\mbox{Diff}\,^m_b$, $\mbox{Diff}\,^m_{sc}$ and $\mbox{Diff}\,^m$ coincide on $\Sigma$).

Analogously, we may define:

\ 

$\mbox{Diff}\,^m_{sc}(X):=$ space of linear operators $P:C^\infty(X) \longrightarrow C^\infty(X)$ such that, locally and with $p_{i,I}\in C^\infty(X)$,
$$P=\Sigma_{i+|I|\leq m}\,p_{i,I}(x,y)(x^2 D_x)^i(x D_y)^I$$
As above, we may also define $\mbox{Diff}\,^m_{sc}(E;F)$, $\mbox{Diff}\,^m_{sc}(E)$. 

Notice that, when one commutes operators, lower-order terms appear; eg, $x(x\partial_x)=(x\partial_x-I)\,x$. However, it is easy to prove the following:

\begin{lemma} $x^\delta\cdot\mbox{Diff}\,^m_b=\mbox{Diff}\,^m_b \cdot x^\delta; x^\delta \cdot \mbox{Diff}\,^m_{sc}=\mbox{Diff}\,^m_{sc} \cdot x^\delta$.

Also: $x^m\cdot\mbox{Diff}^m_b\subseteq \mbox{Diff}^m_{sc}$ (the opposite inclusion is, in general, false).
\end{lemma}

Consider, for example, the exterior derivative operator $d$ acting on $i$-forms. Clearly, $d\in \mbox{Diff}\,^1(\Lambda^i\Sigma;\Lambda^{i+1}\Sigma)$ and $d\in \mbox{Diff}\,^1(\Lambda^i (X);\Lambda^{i+1}(X))$. 

It is easy to show that also the following are true:

\begin{itemize}
\item $d\in \mbox{Diff}\,^1_b(^b\Lambda^i(X);\mbox{$^{b}\Lambda^{i+1}(X)$})$
\item $d\in \mbox{Diff}\,^1_{sc}(\mbox{$^{sc}\Lambda^i(X)$};\mbox{$^{sc}\Lambda^{i+1}(X)$})$
\item $d\in x\cdot\mbox{Diff}\,^1_b(\mbox{$^{sc}\Lambda^i(X)$};\mbox{$^{sc}\Lambda^{i+1}(X)$})$
\end{itemize}
In particular, this shows that there exist $\tilde{d},\overline{d} \in \mbox{Diff}\,^1_b(^{sc}\Lambda^iX;\mbox{$^{sc}\Lambda^{i+1}X$})$ such that $d=x\cdot\tilde{d}=\overline{d}\cdot x$.

We now move on to define some particularly interesting categories of metrics on $\Sigma$.

\begin{definition} An ``exact b-metric'' on $X$ is any metric $g_b$ on $\Sigma$ such that, for some boundary-defining function $x$, $g_b=\frac{d\,x^2}{x^2}+h$, where $h$ denotes a symmetric 2-tensor on $TX$ such that $h_{|\partial X}$ is positive-definite on $T(\partial X)$.

The pair $(X,g_b)$ is an ``asymptotically cylindrical manifold with link $(\partial X,h_{|\partial X})$''.

A ``scattering metric'' on $X$ is any metric $g_{sc}$ on $\Sigma$ such that, for some boundary-defining function $x$, $g_{sc}=\frac{d\,x^2}{x^4}+\frac{h}{x^2}$ where $h$ has the same properties as above.

The pair $(X,g_{sc})$ is a ``scattering manifold with link $(\partial X,h_{|\partial X})$''.
\end{definition}

More explicitly, an exact b-metric has the form
$$g=\frac{d\,x^2}{x^2}+a_{00}d\,x^2+a_{0j}d\,x\,d\,y^j+a_{j0}d\,y^j\,d\,x+a_{ij}d\,y^i\,d\,y^j$$
where $a_{00},a_{0j}=a_{j0},a_{ij}=a_{ji}\in C^\infty(X)$ and $a_{ij}$ has the following property (with respect to its Taylor expansion in $x=0$):
$$a_{ij}(x,y)=g_{ij}(y)+x\,\tilde{a}_{ij}(x,y) \mbox{\ \ where $g_{ij}\,d\,y^id\,y^j$ is a metric on $\partial X$}$$
Notice that, if $g=\frac{dx^2}{x^2}(1+xa)+h$ for some $a\in C^\infty(X)$, the change of variables $\xi:=x+\frac{\gamma}{2}x^2$, where $\gamma:=a(0,y)$, gives $g=\frac{d\xi^2}{\xi^2}+h'$; i.e., such a $g$ is an exact b-metric.

The analogous fact does not hold for sc-metrics. Scattering metrics may thus be generalized as follows.

\begin{definition} An ``asymptotically conical metric'' on $X$ is any metric $g_{ac}$ on $\Sigma$ such that, for some boundary-defining function $x$ and some $a\in C^\infty(X)$, $g_{ac}=\frac{dx^2}{x^4}(1+xa)+\frac{h}{x^2}$, where $h$ denotes a symmetric 2-tensor on $TX$ such that $h_{|\partial X}$ is positive-definite on $T(\partial X)$.

The pair $(X,g_{ac})$ is an ``asymptotically conical manifold with link $(\partial X,h_{|\partial X})$''.
\end{definition}

Clearly, the class of ac-metrics contains the class of sc-metrics as a subset; notice also that, if $g$ is an ac-metric with respect to $x$, $x^2g$ is an exact b-metric.

Fixing one of these three types of metrics implies restricting the choice of boundary-defining functions to those that allow the metric to be brought to the standard form seen in the definition.

Notice that $\nu_b=\{Z\in\Lambda^0(TX):\|Z\|_{g_b}$ is uniformly bounded on $\Sigma\}$.

Analogously,  $\nu_{sc}=\{Z\in\Lambda^0(TX):\|Z\|_{g_{ac}}$ is uniformly bounded on $\Sigma\}$.

\begin{lemma} \label{lemma_local} Let $(X,g)$ be an asymptotically conical manifold, with curvature tensor $R$ and Levi-Civita connection $\nabla$. Then:
\begin{enumerate} 
\item With respect to the basis $\{\partial_x,\partial_{y^i}\}$,
\begin{itemize}
\item $g\in\left( \begin{array}{cc}x^{-4}+C^\infty x^{-3} & C^\infty x^{-2}\\
C^\infty x^{-2} & h_{ij|\partial X}\,x^{-2}+C^\infty x^{-1} \end{array} \right)$ 

where $h_{ij}=h(\partial_{y^i},\partial_{y^j})$.
\item $g^{-1}\in\left(\begin{array}{cc} x^4+C^\infty x^5 & C^\infty\,x^4\\
C^\infty\,x^4 & h^{ij}_{|\partial X}\,x^2+C^\infty\,x^3 \end{array}\right)$

where $h^{ij}$ denotes the inverse matrix of $h_{ij}$ and the notation ``$f\in C^\infty x^\alpha$'' means: $f=\phi\cdot x^\alpha$, for some $\phi\in C^\infty(X)$.
\item The Christoffel symbols have the following property: $\Gamma_{ij}^k\in C^\infty;\Gamma_{ij}^x\in C^\infty\, x;\Gamma_{ix}^k\in C^\infty x^{-1};\Gamma_{ix}^x \in C^\infty\,x;\Gamma_{xx}^k \in C^\infty x^{-1};\Gamma_{xx}^x\in C^\infty x^{-1}$.
\item Let $\tilde{\Gamma}^\gamma_{\alpha\beta}$ denote the Christoffel symbols with respect to the dual basis $\{dx,dy^i\}$; i.e., $\nabla_{\partial_x}dx=\tilde{\Gamma}^x_{xx}dx+\tilde{\Gamma}^i_{xx}dy^i$, etc. Then
$$\tilde{\Gamma}^k_{ij}=(\nabla_{\partial_{y^i}}dy^j)(\partial_{y^k})=\partial_{y^i}[dy^j(\partial_{y^k})]-dy^j(\nabla_{\partial_{y^i}}\partial_{y^k})=-\Gamma^j_{ik}$$
Likewise, $\tilde{\Gamma}^x_{ij}=-\Gamma^j_{ix},\tilde{\Gamma}^k_{xj}=-\Gamma^j_{xk}$, etc.
\end{itemize}
\item $\|R\|\in C^\infty\,x^2$.
\item Let $X\in \nu_b$. Then $\nabla_X\in \mbox{Diff}\,^1_b(^{sc}TX)$ and $\nabla_X\in \mbox{Diff}\,^1_b(^{sc}T^*X)$.

Let $X\in \nu_{sc}$. Then $\nabla_X\in x\cdot \mbox{Diff}\,^1_b(^{sc}TX)$ and $\nabla_X\in x\cdot\mbox{Diff}\,^1_b(^{sc}T^*X)$.
\end{enumerate}
\end{lemma}

Before moving onto the next section, it is probably worth-while making a few final remarks. 

Recall the following

\begin{definition}
Let $\Sigma^n$ be a connected, oriented manifold and $M^{n-1}$ be a (possibly not connected) compact oriented manifold. $\Sigma$ is a ``manifold with ends (with link $M$)'' iff $\Sigma$ admits a decomposition $\Sigma=\Sigma_0 \bigcup_M \Sigma_\infty$, where

\begin{itemize}
\item $\Sigma_0$ is a compact manifold with boundary $\partial(\Sigma_0)=M$
\item $\Sigma_\infty$ is diffeomorphic to $M\times [a,b)$
\end{itemize}
and we are identifying $\partial(\Sigma_0)$ with $M\times\{a\}$.
\end{definition}
\textbf{Remarks:}
\begin{enumerate}
\item Let $\Sigma$ be a manifold with ends. Then $X:=\Sigma_0\bigcup_M M\times[a,b]$ is a compact manifold with boundary and $\Sigma=X\setminus \partial X$. Viceversa, given any compact manifold $X$ with boundary, the local trivialization of $X$ near $\partial X$ given by a boundary-defining function shows that $\Sigma:=X\setminus \partial X$ is a manifold with ends, with link $M:=\partial X$.
\item Notice that any two intervals $[a,b),[a',b')$ are diffeomorphic. The ends of $\Sigma$ can thus be parametrized in countless ways. 

However, not all these diffeomorphisms extend smoothly up to the boundary. Different parametrizations of $\Sigma_\infty$ may thus lead to different differentiable structures on the compactification $X$ of $\Sigma$. When starting out with a manifold with ends, it is important to specify which compactification is being used.
\item Let $(X, g_{ac})$ be asymptotically conical and let $\Sigma$ be the corresponding manifold with ends. The diffeomorphism $x \in (0,1] \simeq r:=\frac{1}{x}\in [1,\infty)$ induces a coordinate system on the ends of $\Sigma$: $\Sigma_\infty \simeq M\times [1,\infty)$.

Notice that $d\,r=-\frac{d\,x}{x^2},\partial_r=-x^2\partial_x$ so that
$$\frac{d\,x^2}{x^4}(1+xa)+\frac{h}{x^2}=d\,r^2+r^2g_M+R$$
where $g_M:=h_{|T(\partial X)}, R:=\frac{a}{x^3}d\,x^2+\frac{h-g_M}{x^2}$. This expression of the metric justifies the name ``asymptotically conical''. In particular,
\begin{eqnarray*}
R(\partial_r,\partial_r)&=&R(x^2\partial_x,x^2\partial_x)=O(r^{-1})\\
R(r^{-1}\partial_{y^i},r^{-1}\partial_{y^j})&=&R(x\partial_{y^i},x\partial_{y^j})=O(r^{-1})\\
R(r^{-1}\partial_{y^i},\partial_r)&=&-R(x\partial_{y^i},x^2\partial_x)=O(r^{-1})
\end{eqnarray*}
so $\|R\|=O(r^{-1})$, where $\|\cdot\|$ is calculated with respect to the metric $d\,r^2+r^2g_M$.

The primary example of scattering manifold is $(\R^n,g_{std})$, compactified via stereographic projection to the half-sphere. See [M2] for details.

It is simple to show that the standard symplectic structure $\omega_{std}$ on $\R^{2n}$ is an element of $^{sc}\Lambda^2 X$.

\item Let $(X,\frac{d\,x^2}{x^2}+h)$ be asymptotically cylindrical and let $\Sigma$ be the corresponding manifold with ends. The diffeomorphism $x\in(0,1]\simeq r:=-log\,x\in [0,\infty)$ induces a coodinate system $\Sigma_\infty\simeq M\times [0,\infty)$ with $d\,r=-\frac{d\,x}{x},\partial_r=-x\,\partial_x$. Thus
$$\frac{d\,x^2}{x^2}+h=d\,r^2+g_M+R$$
where $g_M:=h_{|T(\partial X)}$, $R:=h-g_M$. This expression of the metric justifies the name ``asymptotically cylindrical''.
\item It is interesting to compare these definitions of ``asymptotically conical'' and ``asymptotically cylindrical'' metrics with other definitions available in the literature.

All of them require some form of decay of the ``perturbation term'' $R$, defined in remarks 3,4 above, and of its covariant derivatives; our definitions are however slightly stronger, as they require the tensor $R$ (up to renormalization) to have a smooth extension up to $\partial X$. This allows for a much greater control of the metrics' consequent properties.

Most of the properties relevant to this article would, however, continue to hold for weaker definitions (cfr. eg [MO]).
\end{enumerate}
\section{The Laplace operator on asymptotically conical manifolds} \label{section_laplace}

For any oriented Riemannian manifold $(\Sigma^n,g)$, we will use the following, standard, notation.

\ 

$*_g:\Lambda^i\Sigma \longrightarrow \Lambda^{n-i}\Sigma$ denotes the usual Hodge-star operator.

$d^*_g:\Lambda^i\Sigma \longrightarrow \Lambda^{i-1}\Sigma$ denotes the formal adjoint of $d$. 

$\Lambda^*\Sigma:=\bigoplus \Lambda^i\Sigma$.

\ 

As usual, $d^*_g=(-1)^{n(i+1)+1}*_g d\,*_g$, so $d^*_g\in \mbox{Diff}\,^1(\Lambda^i\Sigma;\Lambda^{i-1}\Sigma)$.

Let $X$ denote an oriented compact manifold with boundary and $\Sigma:=X\setminus \partial X$. For any fixed b-metric on $X$, with volume form $vol_b$, we may define:

\ 

$L^2_b(X):=\{f\in L^2_{loc}(\Sigma) \mbox{ such that } \|f\|_b:=(\int_\Sigma f^2vol_b)^{1/2}<\infty\}$.

$H^m_b(X):=\{f\in L^2_{loc}(\Sigma): P\,f \in L^2_b(X), \  \forall P\in \mbox{Diff}\,^m_b(X)\}\ \ (m\geq 0)$.

$\| f \|_{m,b}:=(\Sigma_{i+|I|\leq m}\| (x\partial_x)^i\partial_y^I\,f\|^2_b)^{1/2}$ defines a norm on $H^m_b(X)$.

\ 

More generally, if $E$ is a metric vector bundle over $X$, we may define:

\ 

$L^2_b(E):=\{f\in L^2_{loc}(E_{|\Sigma}) \mbox{ such that } \int_\Sigma |f|^2 vol_b<\infty\}$. 

$H^m_b(E):=\{f\in L^2_{loc}(E_{|\Sigma}): P\,f \in L^2_b(E), \  \forall P\in \mbox{Diff}\,^m_b(E)\}\ \ (m\geq 0)$.

\ 

Finally, for $\delta \in \R$, we define ``weighted Sobolev spaces'' as follows:

\ 

$x^\delta H^m_b:=\{f\in L^2_{loc}: f=x^\delta u, u\in H^m_b\}$.

$\| f \|_{\delta,m,b}:=\| x^{-\delta}f \|_{m,b}$ defines a norm on $x^\delta H^m_b$, making it isometric to $H^m_b$.

\ 

Analogously, if we endow $X$ with a fixed ac-metric, we may define $L^2_{sc}$, $H^m_{sc}$ and $x^\delta H^m_{sc}$, spaces of functions and sections, using the induced volume form $vol_{ac}$ and operators $P\in \mbox{Diff}\,^m_{sc}$.

Since $\Lambda^i(X)_{|\Sigma}=\mbox{$^b\Lambda^i(X)_{|\Sigma}$}=\mbox{$^{sc}\Lambda^i(X)$}_{|\Sigma}=\Lambda^i(\Sigma)$, when $E$ is one of these bundles one can use the notation $H^k_b(\Lambda^i),H^k_{sc}(\Lambda^i)$. In this case, however, it is important to specify which metric is being used on the bundle.

All the above are Hilbert spaces and contain the space $\Lambda^0_c(E_{|\Sigma})$ as a dense subset. Both ac- and b-metrics have the right properties (bounded curvature and positive injectivity radius) for the standard Sobolev immersion theorems to hold: cfr. [A]. For example, if $f\in x^\delta H^k_b(X)\ \ (k>\frac{n}{2})$, then $x^{-\delta}f$ is continuous and bounded, so $f=O(x^\delta)$.

Notice that, for any constant $c\in \R$, 
$$c\in x^\delta L^2_b(X) \Leftrightarrow x^{-\delta}\in L^2_b\Leftrightarrow \int^1_0 x^{-2\delta-1}d\,x<\infty \Leftrightarrow \delta<0$$
Notice also that, if $g_b=x^2g_{ac}$, then $vol_b=x^n vol_{ac}$, so that $L^2_{sc}=x^{\frac{n}{2}}L^2_b$.

It is a simple consequence of the definitions that any $P\in\mbox{Diff}\,^m_b$ has a continuous extension $P:H^{k+m}_b \longrightarrow H^k_b$. The analogous fact holds for $P\in \mbox{Diff}\,^m_{sc}$. The fact that $\mbox{Diff}\,^m_bx^\delta=x^\delta\mbox{Diff}\,^m_b$ shows however that, for b-metrics (and analogously for ac-metrics), a continuous extension exists also between weighted Sobolev spaces:

\begin{lemma} Let $P\in \mbox{Diff}\,^m_b(E)$. Then, $\forall \delta\in \R, \forall k \geq 0$, $P$ has a continuous extension $P:x^\delta H^{k+m}_b(E) \longrightarrow x^\delta H^k_b(E)$. 

Analogously, let $P\in \mbox{Diff}\,^m_{sc}(E)$. Then $\forall \delta\in \R, \forall k \geq 0$, $P$ has a continuous extension $P:x^\delta H^{k+m}_{sc}(E) \longrightarrow x^\delta H^k_{sc}(E)$.
\end{lemma}

\textbf{Proof}: Consider, for example, $P\in\mbox{Diff}\,^m_b$. Then
$$\|P\,f\|_{\delta,k,b}=\|x^{-\delta}P\,f\|_{k,b}=\|\tilde{P}(x^{-\delta}f)\|_{k,b}\leq C\|x^{-\delta}f\|_{k+m,b}=C\|f\|_{\delta,k+m,b}$$
for some $\tilde{P}\in \mbox{Diff}\,^m_b$, $C>0$.
\qed
Analogously to the compact case, there is a well-developed theory of b- and sc-elliptic operators. We will only need the former.
\begin{definition}Let $P\in \mbox{Diff}\,^m_b(X)$, $P=\Sigma_{i+|I|\leq m}\,p_{i,I}(x,y)(x\,D_x)^i D_y^I$. 

$P$ is b-elliptic iff 
$$\sigma_m(P)(\xi,\eta):=\Sigma_{i+|I|=m}p_{i,I}(x,y)\xi^i\eta^I\neq 0, \ \ \forall (x,y)\in X,\forall (\xi,\eta)\neq 0$$
\end{definition}
It turns out that the properties of a b-elliptic operator $P$ are strongly related to $P_{|\partial X}$. We need the following definition:
\begin{definition} Let $P\in \mbox{Diff}\,^m_b(X)$ be b-elliptic. 

For $\lambda \in \C$, let $\hat{P}(\lambda):C^\infty(\partial X)\longrightarrow C^\infty(\partial X)$ denote the operator locally given by $\hat{P}(\lambda):=\Sigma_{i+|I|\leq m}p_{i,I}(0,y)\lambda^i D_y^I$. This is called the ``indicial operator'' associated to $P$.

Let $spec(P):=\{\lambda\in\C:\hat{P}(\lambda)\mbox{ is not invertible on }C^\infty(\partial X)\}$.
\end{definition}
The above definitions generalize to $P\in \mbox{Diff}\,^m_b(E;F)$. 

In [M1], Melrose constructs a class of pseudo-differential operators which lead to the following result:

\begin{theorem} \label{thm_Fred} Let $P\in \mbox{Diff}\,^m_b(E;F)$ be b-elliptic. Then
\begin{enumerate}
\item $P:x^\delta H^{k+m}_b(E) \longrightarrow x^\delta H^k_b(F)$ is Fredholm iff $\delta \notin -\mathcal{I}m\,spec(P)$
\item $u \in x^\delta H^{k+m}_b(E), \ Pu \in x^\delta H^{k+1}_b(F) \Rightarrow u \in x^\delta H^{k+m+1}_b(E)$
\end{enumerate}
\end{theorem}

The set $-\mathcal{I}m\,spec(P) \subseteq \R$ turns out to be discrete. We call these the ``exceptional weights (of $P$)''. Given a non-exceptional weight $\delta \in \R$, $[\delta]_P$ will denote the connected component of $\R\setminus -\mathcal{I}m \,spec(P)$ containing $\delta$. 

\ 

Having laid out the relevant foundations, we may now focus on the specific operators that will be important further on in this paper.

Let us thus fix an asymptotically conical manifold $(\Sigma^n,g)$. Let $(X,g)$ be its ``scattering compactification''.

Since $^{sc}\Lambda^*X$ is generated by forms of length 1, the definition of $*_g$ shows that $*_g(^{sc}\Lambda^iX)\subseteq \mbox{$^{sc}\Lambda^{n-i}X$}$. Thus $d^*_g$ restricts to an operator in $\mbox{Diff}\,^1_{sc}(^{sc}\Lambda^iX;\mbox{$^{sc}\Lambda^{i-1}X$})$. Actually, if we set $\tilde{\delta}:=(-1)^{n(i+1)+1}*_g\tilde{d}\,*_g$ and $\overline{\delta}:=(-1)^{n(i+1)+1}*_g\overline{d}\,*_g$, we find that 
$$d^*_g=x\cdot\tilde{\delta}=\overline{\delta}\cdot x \in x\cdot \mbox{Diff}\,^1_b(^{sc}\Lambda^iX;\mbox{$^{sc}\Lambda^{i-1}X)$}\subseteq \mbox{Diff}\,^1_{sc}$$
We now define: 

\ 

$D_g:=d\oplus d^*_g\in \mbox{Diff}\,^1_{sc}(^{sc}\Lambda^*X;\mbox{$^{sc}\Lambda^*X$})$.

$\Delta_g:=D_g \circ D_g= d\,d^*_g+d^*_gd \in \mbox{Diff}\,^2_{sc}(^{sc}\Lambda^iX;\mbox{$^{sc}\Lambda^iX$})$.

\ 

Clearly, $D_g,\Delta_g$ are the restrictions to $^{sc}\Lambda^*X$ of the usual operators defined on $\Lambda^*\Sigma$.

\begin{lemma} (cfr. [M2]) Let $(\Sigma,g)$ be asymptotically conical. Then
\begin{enumerate}
\item $\Delta_g=x\cdot(\tilde{d}\oplus \tilde{\delta})\circ (\overline{d}\oplus \overline{\delta})\cdot x\in x^2\cdot\mbox{Diff}\,^2_b(\mbox{$^{sc}\Lambda^iX$};\mbox{$^{sc}\Lambda^iX$})$.
\item On functions, one finds the following expressions:
\begin{itemize}
\item $\Delta_g=x^2\Delta_M+\sqrt{-1}x(n-1)x^2D_x+(x^2D_x)(x^2D_x)+x^3\mbox{Diff}\,^2_b$
\item$(\tilde{d}\oplus \tilde{\delta})\circ (\overline{d}\oplus \overline{\delta})=\Delta_M+(xD_x)(xD_x)+\sqrt{-1}n(xD_x)+1-n+x\mbox{Diff}\,^2_b$
\end{itemize}
\end{enumerate}
\end{lemma}
In particular, the operator $(\tilde{d}\oplus \tilde{\delta})\circ(\overline{d}\oplus \overline{\delta})$ is b-elliptic: for example, on $C^\infty(X)$, the previous lemma shows that 
$$\sigma_2((\tilde{d}\oplus \tilde{\delta})\circ(\overline{d}\oplus \overline{\delta}))(\xi,\eta)=\xi^2+g^{ij}_M\eta^i\eta^j=\xi^2+|\eta|^2$$

We now introduce one last piece of notation: given an asymptotically conical manifold $(\Sigma,g)$, we let $H^k_b(\Lambda^i,sc)$ denote the Sobolev spaces defined above, with respect to the following choices:
\begin{itemize}
\item on $\Lambda^i\Sigma$ we use the metric induced by the ac-metric g;
\item the volume form on $\Sigma$ is the one given by the b-metric $x^2g$.
\end{itemize}
The following result shows that these are the correct spaces in which to study $\Delta_g$, $D_g$.
\begin{corollary} Let $(\Sigma,g)$ be an asymptotically conical manifold. Then

\begin{enumerate} 
\item $\Delta_g:x^\delta H^{k+2}_b(\Lambda^i,sc) \longrightarrow x^{\delta+2} H^k_b(\Lambda^i,sc)\ \ (k \geq 0)$ 

is Fredholm, except for a discrete set of ``exceptional weights''.
\item $D_g:x^\delta H^{k+1}_b(\Lambda^*,sc) \longrightarrow x^{\delta+1} H^k_b(\Lambda^*,sc)\ \ (k\geq 1)$ 

is Fredholm, except for a discrete set of ``exceptional weights''.
\item All elements in $Ker(\Delta_g), Ker(D_g)$ are smooth.
\end{enumerate}
\end{corollary}

\textbf{Proof}: Consider $\Delta_g=x\cdot(\tilde{d}\oplus\tilde{\delta})\circ(\overline{d}\oplus \overline{\delta})\cdot x$\,: 
$$x^\delta H^{k+2}_b(\Lambda^*,sc)\longrightarrow x^{\delta+1}H^{k+2}_b(\Lambda^*,sc)\longrightarrow x^{\delta+1}H^k_b(\Lambda^*,sc) \longrightarrow x^{\delta+2}H^k_b(\Lambda^*,sc)$$
Since multiplication by $x$ is an isometry, by the previous theorem we see that $\Delta_g$ is Fredholm, except for a discrete set of weights. 

Since $D_g=x(\tilde{d}\oplus\tilde{\delta})=(\overline{d}\oplus \overline{\delta})x$, this shows that $D_g:x^{\delta}H^{k+2}_b \longrightarrow x^{\delta+1}H^{k+1}_b$ satisfies dim $Ker(D_g)<\infty$ and that $D_g:x^{\delta+1}H^{k+1}_b \longrightarrow x^{\delta+2}H^{k}_b$ has closed image and finite-dimensional cokernel. Thus $D_g$ is Fredholm if $\delta, \delta-1$ are non-exceptional for $\Delta_g$. 

Finally, to prove (3) notice that, by the above theorem, $Ker(\Delta_g) \subseteq \bigcap_{m\geq 2} x^\delta H^m_b$. Applying the standard Sobolev immersion theorems to $H^m_b$, we get smoothness for the elements of $Ker(\Delta_g)$; (3) then follows from $Ker(D_g)\subseteq Ker(\Delta_g)$. \qed
We may thus define:

\ 

$\mathcal{H}^i_\delta(\Sigma):=Ker(\Delta_g)$, where $\Delta_g:x^\delta H^{k+2}_b(\Lambda^i,sc) \longrightarrow x^{\delta+2} H^k_b(\Lambda^i,sc)$

$\mathcal{K}^i_\delta(\Sigma):=Ker(D_g)$, where $D_g:x^\delta H^{k+1}_b(\Lambda^i,sc) \longrightarrow x^{\delta+1} H^k_b(\Lambda^*,sc)$

\ 

From the corollary, we see that these spaces are independent of $k$. When $\Delta_g$ and $D_g$ are Fredholm, they have Fredholm indices $i_{\Delta_{g,\delta}}$ and $i_{D_{g,\delta}}$. The following lemma investigates the dependence of these spaces and indices on $\delta$.
\begin{lemma} Let $(\Sigma,g)$ be an asymptotically conical manifold. Then
\begin{enumerate}
\item $i_{\Delta_{g,\delta}}$ and $\mathcal{H}^i_\delta(\Sigma)$ depend only on $[\delta]_\Delta$.
\item $i_{D_{g,\delta}}$ and $\mathcal{K}^i_\delta(\Sigma)$ depend only on $[\delta]_D$.
\end{enumerate}
\end{lemma}
The proof is based on the invariance of the Fredholm index of curves of Fredholm operators.
\subsection{$\Delta_g$ on functions} \label{subsection_functions}

When studying $\Delta_g$ on functions, one gets extra information thanks to the elliptic maximum principle.

\begin{prop} Let $(\Sigma,g)$ be an asymptotically conical manifold.

Consider $\Delta_g:x^\delta H^{k+2}_b(\Sigma) \longrightarrow x^{\delta+2}H^k_b(\Sigma)$, for $\delta$ non-exceptional. Then:
\begin{enumerate}
\item $\delta>0 \Rightarrow \mathcal{H}^0_\delta=0$. Thus:
\begin{itemize}
\item $\Delta_g$ is injective.
\item $d: x^\delta H^{k+2}_b(\Sigma) \longrightarrow x^{\delta+1} H^{k+1}_b(\Lambda^1,sc)$ is injective.
\end{itemize}
\item $\delta<n-2 \Rightarrow Coker(\Delta_g)=0$. Thus:
\begin{itemize}
\item $\Delta_g$ is surjective.
\item $d^*_g: x^\delta H^{k+1}_b(\Lambda^1,sc) \longrightarrow x^{\delta+1} H^k_b(\Sigma)$ is surjective for $\delta<n-1$ non-exceptional.
\end{itemize}
\item Dim $\mathcal{H}^0_\delta(\Sigma)$ is independent of the particular choice of asymptotically conical metric on $\Sigma$ (with respect to a fixed link $(M,g_M)$).
\end{enumerate}
\end{prop}
The proof is based on the elliptic maximum principle and the self-adjointness of $\Delta_g$ on $C^{\infty}_c(\Sigma)$. Cfr. [CZ] for similar statements and techniques. Notice that, in particular, $\Delta_g$ is an isomorphism for $0<\delta<n-2$.

The above proposition has an interesting consequence:

\begin{lemma}(``gluing principle for harmonic functions'')
 Let $(\Sigma,g)$ be asymptotically conical. For $\delta < n-2$, let $f$ be a smooth function in $x^\delta H^k_b(\Sigma)$ such that $\Delta_g f_{|\Sigma^i_\infty} = 0$. 

Then there exists a unique $F \in x^\delta H^k_b(\Sigma)$ such that $|F(x)-f(x)|\rightarrow 0$ (as $x\rightarrow 0$) and $\Delta_g F \equiv 0$.
\end{lemma}

\textbf{Proof}: Since $\Delta_g f$ is smooth and has compact support, it is clear that, $\forall \eta \in \R, \forall s\geq 0,\ \ \Delta_g f \in x^{\eta+2} H^s_b(\Sigma)$.

Fix any $\eta: \mbox{max}\{0,\delta\}<\eta<n-2$. Since $\Delta_g:x^\eta H^{s+2}_b\longrightarrow x^{\eta+2} H^s_b$ is surjective, there exists $\tilde{f} \in x^\eta H^{s+2}_b: \Delta_g \tilde{f}=\Delta_g f$. In particular, $\tilde{f} \in x^\delta H^{s+2}_b$. Choosing $s$ large enough, we get $|\tilde{f}|=O(x^\eta)$, so $\tilde{f}\rightarrow 0$.

It is now enough to set $F:=f-\tilde{f}$. Uniqueness is a consequence of the elliptic maximum principle.
\qed
In other words, any collection of harmonic functions (eg: constants) on $\Sigma^i_\infty$ may, up to a slight perturbation, be ``glued together'' to get a harmonic function on $\Sigma$ with the same asymptotic behaviour.

Going a step further, [M1] provides a complete asymptotic expansion of the harmonic functions on any asymptotically conical $(\Sigma,g)$ in terms of the metric on the link. In particular, Christiansen and Zworski ([CZ]) show how this can be used to relate the harmonic functions on ($\Sigma,g)$ to the eigenfunctions on the link. We may summarize their results as follows.

\begin{itemize}
\item Recall that, on $(M^{n-1}\times (1,\infty), g:=d\,r^2+r^2\,g_M)$, 
$$\Delta_g=\frac{\Delta_{g_M}}{r^2}-\frac{n-1}{r}\partial_r-\partial_r^2$$

Let $f\in C^\infty(M^{n-1})$ be any eigenfunction of $\Delta_{g_M}$, relative to any eigenvalue $\lambda_i$. Then it is easy to check that the function $f\,r^{a_i}$ is harmonic on $M\times (1,\infty)$, where $a_i:=\frac{2-n+\sqrt{(2-n)^2+4\lambda_i}}{2}$.

Notice that, when $(M,g_M)=(S^{n-1},g_{std})$, these harmonic functions on $\R^n\setminus B$ are exactly the homogeneous harmonic polynomials.
\item If $(\Sigma,g)$ has link $(M,g_M)$, these functions are ``asymptotic models'' for the harmonic functions (with polynomial growth) on $\Sigma$, in the following sense: given any $f\in \mathcal{H}^0_\delta(\Sigma)$, on each end $f$ converges to a linear combination of the above.
\item Viceversa, assigning an ``asymptotic model'' to each end determines a unique $f\in \mathcal{H}^0_\delta$ with that asymptotic behaviour.
\end{itemize}

This allows us to express dim $\mathcal{H}^0_\delta$ in terms of the number of ends and the dimension of the space of eigenfunctions on the components $(M_i,g_{M_i})$ of the link $(M,g_M)$.

In particular, we will be interested in the following conclusion.

\begin{definition} We say that a harmonic function $f$ has ``strictly sub-linear growth'' if it is asymptotic to a linear combination of models $f_ir^{a_i}$ with $0\leq a_i<1$ and at least one $a_i>0$. In other words, if there exist $\delta>-1, \epsilon>0$ such that $f\in \mathcal{H}^0_\delta(\Sigma)$, $f\notin\mathcal{H}^0_{-\epsilon}(\Sigma)$.
\end{definition}

\begin{corollary} Let $(\Sigma,g)$ be asymptotically conical.

Then there exist harmonic functions with strictly sub-linear growth iff at least one end has an eigenvalue in the interval $(0,n-1)$.
\end{corollary}

\textbf{Proof}: Notice that $0<a_i<1 \Leftrightarrow 0<\lambda_i<n-1$. We may use any asymptotic model of the type $\Sigma_{i\neq j}c_j+f_ir^{a_i}$, where $c_j\in\R$.
\qed

\subsection{$D_g$ on $1$-forms} \label{subsection_1forms}

We are now going to calculate dim $\mathcal{K}^1_\delta=$dim $Ker(D_g)$, where $D_g$ is acting on weighted 1-forms, in two cases of particular interest.

In this section, $i$ will denote the canonical map $i:H^1_c(\Sigma) \longrightarrow H^1(\Sigma)$.
 
We start with the following

\begin{lemma} Let $(\Sigma,g)$ be asymptotically conical.

For any $\epsilon>0$, there exists an injective map
$$q:H^1(\Sigma) \longrightarrow \mathcal{K}^1_{1-\epsilon}$$
which is independent of $\epsilon$. With respect to any $\delta\in [1,n-1)$, it has the following property:
$$q([\alpha])\in \mathcal{K}^1_\delta \Leftrightarrow [\alpha]\in i(H^1_c)$$
\end{lemma}
\textbf{Proof}: Begin by considering any $[\alpha] \in H^1(\Sigma)$. We first show that it has a representative in $\mathcal{K}^1_{1-\epsilon}$.

Recall the isomorphism
$$\begin{array}{rcl}
H^1(M)&\simeq&H^1(M\times[1,\infty))\\
\mbox{$[\beta_M]$} &\mapsto& \mbox{$[\beta_\infty]$}
\end{array}$$
where $\beta_\infty$ is defined by $\beta_\infty[\omega,r]:=\beta_M[\omega]$. We call such a form ``translation-invariant''.

Applying this to $[\alpha_{|\Sigma_\infty}]\in H^1(M\times[1,\infty))$ gives a $1$-form $\alpha_M$ on $M$ and a translation-invariant $1$-form $\alpha_\infty$ on $M\times[1,\infty)$ such that $[\alpha_\infty]=[\alpha_{|\Sigma_\infty}]$. In other words,
$$\alpha_{|\Sigma_\infty}=\alpha_\infty+d\,(f_\infty)$$
for some $f_\infty\in C^\infty(\Sigma_\infty)$ and some translation-invariant $\alpha_\infty$.

We may now extend $f_\infty$ to a function $f \in C^\infty(\Sigma)$; this gives an extension of $\alpha_\infty$ to a global $1$-form $\tilde{\alpha}$ on $\Sigma$, defined as $\tilde{\alpha}:=\alpha-d\,f$.

The form $\tilde{\alpha}$ has the following properties:
\begin{itemize}
\item $d\,\tilde{\alpha}=0$: this is clear.
\item $[\tilde{\alpha}]=[\alpha]$: this is clear.
\item $\tilde{\alpha} \in x^\eta H^1_b(\Lambda^1,sc)$, $\forall \eta<1$: locally near $\partial X$, $\tilde{\alpha}=b_i(y)d\,y^i=xb_i(y)\frac{d\,y^i}{x}$. Since $xb_i(y)\in C^\infty(X)$, this shows that $\tilde{\alpha}\in\mbox{$^{sc}\Lambda^1(X)$}$. 

By definition, $\tilde{\alpha}\in x^\eta H^1_b \Leftrightarrow \|x^{-\eta}\tilde{\alpha}\|_{1,b}<\infty$. Recall:
$$\|x^{-\eta}\tilde{\alpha}\|_{1,b}=\|x^{-\eta}\tilde{\alpha}\|_b+\Sigma_i\|x\partial_x(x^{-\eta+1}b_i(y))\frac{d\,y^i}{x}\|_b+\Sigma_{ij}\|\partial_{y^j}(x^{-\eta+1}b_i(y))\frac{d\,y^i}{x}\|_b$$
It is now enough to examine these terms one by one, to get the result. For example:
\begin{eqnarray*}
\|x^{-\eta}\tilde{\alpha}\|^2_b&=&\int_\Sigma\|x^{-\eta}\tilde{\alpha}\|^2_{ac}vol_b\leq c_1+c_2\int^1_0\|x^{-\eta+1}\frac{d\,y^i}{x}\|^2_{ac}\,x^{-1}d\,x\\
&\leq& c_1+c_3\int^1_0x^{-2\eta+1}d\,x
\end{eqnarray*}
so $\|x^{-\eta}\tilde{\alpha}\|_b<\infty$ for $\eta<1$.
\item For $\delta\geq 1$, $\tilde{\alpha}\in x^\delta H^1_b(\Lambda^1,sc)\Leftrightarrow \tilde{\alpha}\in\Lambda^1_c(\Sigma)$: this should be clear from the above.
\end{itemize}

\ 

Consider now, for $\eta<n-1$, the sequence:
$$\CD
x^{\eta-1}H^2_b(\Sigma) @>d>> x^\eta H^1_b @>d^*_g>> x^{\eta+1} L^2_b(\Sigma)
\endCD$$
By surjectivity, setting $\eta=1-\epsilon$, there exists a function $k\in x^{-\epsilon}H^2_b:\Delta_g k=d^*_g\tilde{\alpha}$. Since $\Delta$ is not injective, the choice of $k$ is not canonical; however, $k$ is unique up to Ker($\Delta$) and this is independent of $\epsilon$, so $k$ also is. Notice now that $\tilde{\alpha}-d\,k\in\mathcal{K}^1_{1-\epsilon}$ and that $[\tilde{\alpha}-d\,k]=[\alpha]$. When $\tilde{\alpha}\in \Lambda^1_c(\Sigma)$, we can act in a similar way, choosing however $k$ with respect to $\eta\in[1,n-1)$. 

The above process thus shows how to build the required linear map $q$, defined by $q([\alpha]):=\tilde{\alpha}-d\,k$. To show that $q$ is injective, let $[\alpha_1],...,[\alpha_p]$ be a basis of $H^1(\Sigma)$. Assume $\exists \lambda_i \in \R: \sum\lambda_i(\tilde{\alpha}_i-d\,k_i)=0$. Then $$\sum\lambda_i[\alpha_i]=\sum\lambda_i[\tilde{\alpha}_i-d\,k_i]=0$$
so $\lambda_i=0$. 
\qed

\ 

For all $\eta \in \R$, let $E_\eta:=\{$exact $1$-forms in $\mathcal{K}^1_\eta\}$.

\begin{lemma}
\ 
\begin{enumerate}
\item For $\eta<1$, $E_\eta=d(\mathcal{H}^0_{\eta-1})$.

In particular, dim $E_\eta=$dim $\mathcal{H}^0_{\eta-1}-1$.
\item For $\epsilon>0$ (sufficiently small) and $1-\epsilon<\delta<n-1$, $E_\delta=d(\mathcal{H}^0_{-\epsilon})$.

In particular, dim $E_\delta=s-1$, where $s$ is the number of ends of $\Sigma$.

\end{enumerate}
\end{lemma}
\textbf{Proof}: It is clear that, if $f\in \mathcal{H}^0_{\eta-1}$, then $d\,f \in E_\eta$. Viceversa, assume $d\,f \in E_\eta$. Then $\Delta_g f\equiv 0$, so the asymptotic expansion of $f$ shows that $f\in\mathcal{H}^0_{\eta-1}$. 

The proof of (2) is similar: by the ``gluing principle'', if $f\in\mathcal{H}^0_{-\epsilon}$, then $f=c+\tilde{f}$, where $c$ is a function constant on each end and $\tilde{f}\in x^{n-2-\tilde{\epsilon}}H^s_b$. Thus $d\,f\in E_{n-1-\tilde{\epsilon}}\subseteq E_\delta$. Viceversa, if $d\,f\in E_\delta$, then $f\in d(\mathcal{H}^0_{-\epsilon})$.

Notice that, in calculating dim $E_\eta$ ($\forall \eta\in \R)$, one must take into account the fact that any $c\in \R$ is a harmonic function but $d(c)\equiv 0$, so it does not contribute to $E_\eta$.
\qed

We now have a good picture of the structure of $\mathcal{K}^1_\delta(\Sigma)$:
\begin{prop} \ 
\begin{enumerate}
\item for $\eta<1$, $\mathcal{K}^1_\eta=E_\eta\oplus q(H^1(\Sigma))$
\item for $\delta \in [1,n-1)$, $\mathcal{K}^1_\delta=E_\delta\oplus q(i(H^1_c(\Sigma)))$
\end{enumerate}
\end{prop}
\textbf{Proof}: It is clear from the definitions of $q$ and $E_\eta$ that $E_\eta \bigcap q(H^1(\Sigma))=\{0\}$.
 
Let $H_\eta$ be any subspace of $\mathcal{K}^1_\eta$ containing $q(H^1(\Sigma))$, such that $\mathcal{K}^1_\eta= E_\eta \oplus H_\eta$.

The map $H_\eta\longrightarrow H^1(\Sigma)$, $\alpha \mapsto [\alpha]$ is clearly injective. Thus dim $H_\eta \leq $dim $H^1(\Sigma)$. Since dim $H^1(\Sigma)=$dim $q(H^1(\Sigma))$, we get $H_\eta\leq q(H^1(\Sigma))$, hence the equality.

To prove (2), notice that
$$E_\delta \oplus q(i(H^1_c)) \leq \mathcal{K}^1_\delta \leq E_{1-\epsilon} \oplus q(H^1)$$
Any $\alpha\in \mathcal{K}^1_\delta$ may be written as $\alpha=d\,f+\alpha'$ where $d\,f\in E_{1-\epsilon}, \, \alpha'= q([\alpha]) \in \mathcal{K}^1_{1-\epsilon}$.
Recall that, for small $\epsilon\geq 0$, $E_{1-\epsilon}=E_\delta$. Thus $d\,f\in E_\delta$, which implies that $\alpha'\in\mathcal{K}^1_\delta$ and thus $[\alpha]\in i(H^1_c)$. This proves that $\mathcal{K}^1_\delta\leq E_\delta\oplus q(i(H^1_c))$.
\qed

We will be particularly interested in the following conclusion:

\begin{corollary} Let $(\Sigma,g)$ be asymptotically conical. Then
\begin{enumerate}
\item dim $\mathcal{K}^1_\epsilon=$dim $H^1(\Sigma)+$dim $\mathcal{H}^0_{-1+\epsilon}-1$.

This is the dimension of the space of all $1$-forms $\alpha \in Ker(D_g): \|\alpha\|_g$ decays.
\item $\forall \delta \in (1,n-1)$, dim $\mathcal{K}^1_\delta=$dim $H^1_c(\Sigma)$.

In particular, this shows that the space $Ker(D_{g|L^2_{sc}(\Lambda^1)})=\mathcal{K}^1_{\frac{n}{2}}$ of closed and co-closed 1-forms in $L^2$ has dimension $H^1_c(\Sigma)$.
\end{enumerate}
\end{corollary}
\textbf{Proof}: (1) follows directly from the above.

To prove (2), notice that the previous proposition shows that
$$\mbox{dim }\mathcal{K}^1_\delta=\mbox{dim }i(H^1_c)+s-1$$
where $s$ is the number of ends of $\Sigma$.

Let $M$ be the link of $\Sigma=\Sigma_0\bigcup_M\Sigma_\infty$. The long exact sequence
$$\CD
H^0_c(\Sigma_0) @>>> H^0(\Sigma_0) @>>> H^0(M) @>>> H^1_c(\Sigma_0) @>i>> H^1(\Sigma_0)\\
@|                   @|                 @|          @|                    @|\\
0               @.   \R            @.   \R^s   @.   H^1_c(\Sigma)   @.    H^1(\Sigma)
\endCD$$
now shows that $H^1_c \simeq i(H^1_c) \oplus \R^{s-1}$.
\qed

\section{Deformations of asymptotically conical submanifolds} \label{section_immersions}
Having developed the theory of asymptotically conical manifolds, we now turn to studying submanifolds. The goal of this section is to present a notion of immersion in the category of manifolds with boundary and to study the stability of this condition under deformations. For the purposes of section \ref{section_SLdefs}, the main result is corollary \ref{cor_immersions}, which shows that asymptotically conical submanifolds remain asymptotically conical under ``deformations which keep infinity fixed''.

$X,X'$ will denote compact manifolds with boundary.

$\Sigma:=X\setminus\partial X, N:=X'\setminus\partial X'$ will be their interior.

\begin{definition} An immersion $\phi:X\longrightarrow X'$ is a ``b-immersion''
if
\begin{enumerate}
\item $\phi(\partial X)\subseteq \partial X'$
\item $\forall p \in \partial X$, consider $d\phi[p]:T_pX\longrightarrow T_{\phi(p)}X'$ and, using condition (1), the quotient map 
$$d\phi[p]:T_pX/T_p\partial X \longrightarrow T_{\phi(p)}X'/T_{\phi(p)}\partial X'$$
The condition is that this quotient map be injective.
\end{enumerate}
\end{definition}

Let $\phi:X\longrightarrow X'$ be a b-immersion. Let $x$ be a boundary-defining function for $X'$, $p\in \partial X$ and $V\in T_pX:V\notin T_p\partial X$. Then, by hypothesis,
$$d(x\circ \phi)[p](V)=dx[\phi(p)]d\phi[p](V)\neq 0$$
so $x\circ \phi$ is a boundary-defining function for $X$. From now on, we will often identify $X$ with its image in $X'$, omitting $\phi$ from the notation.

Clearly, $^bTX\leq \mbox{$^bTX'_{|X}$}$ and $^{sc}TX\leq \mbox{$^{sc}TX'_{|X}$}$. 
We may define:

\ 

$\nu_b(X,X'):=\{Z_{|X}:Z\in \nu_b(X')\}$

$\nu_{sc}(X,X'):=\{Z_{|X}:Z\in \nu_{sc}(X')\}$

\

Notice that 

\ 

$\nu_b(X)=\{Z\in \nu_b(X,X'):Z\in \mbox{$^bTX$}\}$

$\nu_{sc}(X)=\{Z\in\nu_{sc}(X,X'):Z\in \mbox{$^{sc}TX$}\}$

\ 

\textbf{Remark:} Analogously to the standard Riemannian case, any b-metric on $X'$ can be seen as a smooth section of a symmetric product of the bundle $^bT^*X'$. It thus restricts to a b-metric on $X$. Regarding sc- and ac-metrics, pointwise there is no difference between them: both come from global sections of the same symmetric product of $^{sc}T^*X'$. The difference arises from the section itself, i.e. from its Taylor expansion at the boundary. Thus, depending on the immersion, a sc-metric on $X'$ may restrict to either a sc- or a ac-metric on $X$.

\ 

Let us fix an asymptotically conical manifold $(X',g)$. From now on, $\nabla, R,exp$ will denote, respectively, the Levi-Civita connection, curvature tensor and exponential map of $(N,g)$. If $\phi$ is a b-immersion and $p\in X$, we will let $\perp, T$ denote the normal and tangential components of any $V\in T_pX'$ with respect to $T_pX$.

We now want to deform $\Sigma\subseteq N$, using the exponential map of $N$. In particular, we now want to prove that, for any $V\in \nu_{sc}(X')$, $exp\,V\circ\phi:X\longrightarrow X'$ is also a b-immersion. The main difficulty lies in the fact that, since $g$ ``blows up'' on $\partial X'$, $exp\,V$ is, a priori, defined only on $N$. We thus need the following

\begin{prop} \label{prop_exp} Let $(X',g)$ be an asymptotically conical manifold. 
\begin{enumerate}
\item Let $V\in \nu_b(X')$. Then $exp\,V:N\longrightarrow N$ has a smooth extension $exp\,V:X'\longrightarrow X'$ such that $exp\,V(\partial X')\subseteq\partial X'$.
\item Let $V\in \nu_{sc}(X')$. Then the extension of $exp\,V$ satisfies $exp\,V_{|\partial X'}\equiv id$ and is an immersion on $\partial X'$.
\item Let $\|\cdot\|_{1,b}$ denote the norm on $H^1_b(TX',sc)$. Then $\exists \delta>0:\forall V\in\nu_{sc}(X'),\|V\|_{1,b}\leq\delta$, $exp\,V:X'\longrightarrow X'$ is an immersion.
\item Let $V\in\nu_{sc}(X')$. Then $exp\,V_*(\nu_{sc})\subseteq \nu_{sc}$.
\end{enumerate}
\end{prop}
\textbf{Proof}: Let $x,y^1,\dots,y^N$ be coordinates on $X'$ and let $V=V^x\,x\partial_x+V^i\,\partial_{y^i}\in \nu_b(X')$. Consider the system of ODE's
$$\left\{ \begin{array}{l}
\ddot{c}^x+c^x(\frac{\Gamma_{ij}^x}{x})(\gamma)\dot{c}^i\dot{c}^j+\Gamma_{xj}^x(\gamma)\dot{c}^x\dot{c}^j+\Gamma_{ix}^x(\gamma)\dot{c}^i\dot{c}^x+\frac{(x\Gamma_{xx}^x)(\gamma)}{c^x}\dot{c}^x\dot{c}^x=0\\
\ddot{c}^k+\Gamma_{ij}^k(\gamma)\dot{c}^i\dot{c}^j+\frac{(x\Gamma_{xj}^k)(\gamma)}{c^x}\dot{c}^x\dot{c}^j+\frac{(x\Gamma_{ix}^k)(\gamma)}{c^x}\dot{c}^i\dot{c}^x+\frac{(x^2\Gamma_{xx}^k)(\gamma)}{(c^x)^2}\dot{c}^x\dot{c}^x=0\\
c(0)=(1,y_0)\\
\dot{c}(0)=(V^x[x_0,y_0],V^1[x_0,y_0],\dots,V^N[x_0,y_0])
\end{array} \right.$$
where:
\begin{itemize}
\item $c^x(s),c^1(s),\dots,c^N(s)$ are the unknown functions;
\item $\Gamma^\sigma_{\mu\nu}$ are the Christoffel symbols of $g$ with respect to $\{\partial_x,\partial_{y^i}\}$;
\item $\gamma(s):=(x_0\,c^x(s),c^1(s),\dots,c^N(s))$;
\item the system depends on the parameters $(x_0,y_0)\in [0,1]\times\partial X'$.
\end{itemize}
By lemma \ref{lemma_local}, all coefficients in the system are smooth (except when $c^x=0$); so, with our initial conditions, there exists a unique solution $c(s)$ and it is smooth with respect to $s,x_0,y_0$. Thus $\gamma(s)$ is well-defined for all $(x_0,y_0)$ and is smooth.

Notice that, for example, $(\frac{\Gamma_{ij}^x}{x})(\gamma)=\frac{\Gamma_{ij}^x(\gamma)}{x_0c^x}$. Thus, multiplying the first equation by $x_0$ ($x_0\neq 0$), the system can be re-written in terms of $\gamma$, showing that $\gamma(s)$ is the geodesic of $(N,g)$ with initial conditions
$$\gamma(0)=(x_0,y_0),\ \ \dot{\gamma}(0)=(x_0V^x[x_0,y_0],V^1[x_0,y_0],\dots,V^N[x_0,y_0])=V[x_0,y_0]$$
This shows that the geodesics generated by any b-vector field, a priori defined only on $N$, extend smoothly, as curves, to $\partial X'$. Since $N$ is complete, it also shows, for $x_0\neq 0$, that $c(s)$ is defined $\forall s\in \R$ and that $c(s)\subseteq N$, since this is true for $\gamma(s)$.

On the other hand, notice that, for initial conditions $(0,y_0)\in \partial X'$, the corresponding curve $\gamma(s)$ is completely contained in $\partial X'$. This proves (1).

Now assume $V\in \nu_{sc}(X')$. Then, for $x_0=0$, $\dot{c}(0)=0$. Thus $c(s)\equiv (1,y_0)$ is the solution, so $\gamma(s)\equiv (0,y_0)$. This proves that, when $V\in\nu_{sc}$, $exp\,V_{|\partial X'}\equiv id$.

Let $p\in \partial X'$. Since $d(exp\,V)[p]_{|T_p\partial X'}\equiv id$, to finish proving (2) it remains only to prove that $d(exp\,V)[p](\partial_x)$ has a component in the $\partial_x$-direction.

In general, let $p\in N$ and $Z\in T_pX'$. Let $c(t)$ be a curve in $N:c(0)=p,\dot{c}(0)=Z$. Let $t\mapsto \gamma(t,s)$ denote the 1-parameter family of geodesics defined by $\gamma(t,0)=c(t),\frac{\partial}{\partial s}\gamma(t,0)=V[c(t)]$.

Let $J(s)$ be the Jacobi vector field solution of 
$$\nabla_s\nabla_s J(s)+R(J(s),\dot{\gamma})\dot{\gamma}=0,\ \ J(0)=Z(p),\ \nabla_sJ(0)=\nabla_Z V(p)$$
where $\dot{\gamma}:=\frac{\partial}{\partial s}\gamma(0,s)$. Notice that
\begin{itemize}
\item $\frac{\partial}{\partial t}\gamma(0,0)=Z(p)$
\item $\nabla_s\frac{\partial}{\partial t}\gamma(0,0)=\nabla_s\nabla_t\gamma(0,0)=\nabla_t\nabla_s\gamma(0,0)=\nabla_t\frac{\partial}{\partial s}\gamma(0,0)=\nabla_ZV(p)$
\item $\nabla_s\nabla_s\frac{\partial}{\partial t}\gamma(0,s)=\nabla_s\nabla_t\frac{\partial}{\partial s}\gamma(0,s)=\nabla_t\nabla_s\frac{\partial}{\partial s}\gamma + R(\frac{\partial}{\partial s}\gamma,\frac{\partial}{\partial t}\gamma)\frac{\partial}{\partial s}\gamma =$

$=-R(\frac{\partial}{\partial t}\gamma(0,s),\dot{\gamma})\dot{\gamma}$
\end{itemize}
This proves that $J(s)\equiv \frac{\partial}{\partial t}\gamma(t,s)_{|t=0}$. Thus, in general, $exp\,V_*$ can be expressed in terms of Jacobi vector fields:
$$exp\,V_*Z=exp\,V_*(\frac{\partial}{\partial t}\gamma(0,0))=\frac{\partial}{\partial t}\gamma (0,1)=J(1)$$
Now, for example, let $Z=x^2\partial_x$. Using the following lemma, we find that
$$J(s)=Z(p)+s\nabla_ZV(p)+Q(s),\ \ \|Q(s)\|=O(x^2)$$
where $Z(p),\nabla_ZV(p)$ are extended along $\gamma$ by parallel transport. Again, we have to prove that this formula has meaning up to $\partial X'$.

Since $exp\,V:X'\longrightarrow X'$ is smooth, $J(1)$ has a smooth extension up to $\partial X'$. Notice also that $Z,\nabla_ZV$ have smooth extensions up to $\partial X'$ (using lemma \ref{lemma_local}).

Recall that the parallel transport of, for example, $Z$ along $\gamma$ is defined by
$$\nabla_sZ(s)=0,\ \ Z(0)=Z[x_0,y_0]$$
Writing $Z=Z^x\partial_x+Z^i\partial_{y^i}$, we find
\begin{eqnarray*}
\nabla_sZ(s) &=& \dot{Z}^x\partial_x+\dot{Z}^i\partial_{y^i}+Z^x\dot{\gamma}^x(\nabla_{\partial_x}\partial_x)(\gamma)+Z^x\dot{\gamma}^k(\nabla_{\partial_{y^k}}\partial_x)(\gamma)\\
&&+Z^i\dot{\gamma}^x(\nabla_{\partial_x}\partial_{y^i})(\gamma)+Z^i\dot{\gamma}^k(\nabla_{\partial_{y^k}}\partial_{y^i})(\gamma)
\end{eqnarray*}
A proof similar to the one above shows that parallel transport extends smoothly up to the boundary, so $Q(1)$ is also smooth up to $\partial X'$.

Writing $Q(1)=a\,\partial_x+b_i\,\partial_{y^i}=\frac{a}{x^2}(x^2\partial_x)+\frac{b_i}{x}(x\partial_{y^i})$, the estimate on $\|Q(1)\|$ shows that $\frac{a}{x^2}=O(x^2)$, $\frac{b_i}{x}=O(x^2)$. A Taylor expansion of $a,b_i$ based at $p\in \partial X'$ thus shows that $Q(1)\in x^2\cdot \nu_{sc}$, so
$$exp\,V_*(Z)=Z+x\cdot\nu_{sc}+x^2\cdot\nu_{sc}$$
Finally:
$$exp\,V_*(\partial_x)=\frac{exp\,V_*(x^2\partial_x)}{x^2}=\frac{x^2\partial_x+x^2\cdot\nu_b+x^2\cdot\nu_{sc}}{x^2}=\partial_x+\nu_b+\nu_{sc}$$
which concludes the proof of (2).

Notice that, in a similar way, one can study the case $Z=x\partial_{y^i}$, showing that $exp\,V_*(\nu_{sc})\subseteq\nu_{sc}$. This proves (4).

To prove (3), we start from the formula
$$exp\,V_*(\partial_x)=\partial_x+\frac{\nabla_{x^2\partial_x}V}{x^2}+\frac{Q(1)}{x^2}$$
From the above, we see that we may write
$$\frac{\nabla_{x^2\partial_x}V}{x^2}=a(x\partial_x)+b_i\partial_{y^i},\ \ \frac{Q(1)}{x^2}=c(x^2\partial_x)+d_i(x \partial_{y^i})$$
Since $\nabla_{x\partial_x}\in \mbox{Diff}\,^1_b$, we find $a\leq \|\frac{\nabla_{x^2\partial_x}V}{x^2}\|_b=\|\nabla_{x\partial_x}V\|_{ac}\leq\|V\|_{1,b}$.

Analogously, one can deduce from lemma \ref{lemma_Jacobi} that $\|Q(1)\|_{ac}\leq x^2\|V\|_{1,b}$, so $c\leq \|\frac{Q(1)}{x^2}\|_{ac}\leq \|V\|_{1,b}$.

In other words, a $\|\cdot\|_{1,b}$-bound yields bounds on $a,c$ which are independent of the particular $V$. This proves that there exists a neighborhood $U$ of $\partial X'$, independent of $V$, such that $exp\,V_*(\partial_x)$ has a non-zero component in the $\partial_x$-direction.

In a similar way, one can study $exp\,V_*(\partial_{y^i})$. The end result is that any $\|\cdot\|_{1,b}$-bound gives a neighborhood $U$ of $\partial X'$, independent of $V$, such that $exp\,V$ is an immersion on $U$.

Since the complement $K$ of $U$ in $X'$ is compact, there exists $\delta>0$ such that $\|V\|_{ac}\leq \delta$ implies that $exp\,V$ is an immersion on $K$. This proves that, for any $V$ such that $\|V\|_{1,b}\leq \delta$, $exp\,V$ is an immersion on $X'$.
\qed
\begin{corollary} \label{cor_immersions} Let $(X',g)$ be an asymptotically conical manifold and let $\phi:X\longrightarrow X'$ be a b-immersion, so that $(X,\phi^*g)$ is asymptotically conical.

Let $\|\cdot\|_{1,b}$ denote the norm on $H^1_b(TX',sc)$. Then $\exists \delta>0:\forall V\in \nu_{sc}(X'), \|V\|_{1,b}\leq \delta$, $\exp\,V\circ \phi$ is a b-immersion and thus also induces an asymptotically conical metric on $X$.
\end{corollary}
\textbf{Remark}: More generally, if we let $\nu_{1,sc}(X')$ denote the $C^1$-analogue of $\nu_{sc}(X')$, the same methods show that
$$V\in\nu_{1,sc}(X'),\|V\|_{1,b}\leq \delta\Rightarrow exp\,V\mbox{ is a $C^1$-immersion up to $\partial X'$}$$
Thus, if $\phi:X\longrightarrow X'$ is a b-immersion, $exp\,V\circ\phi$ is a $C^1$ b-immersion.

\ 

The following lemma, used in the proof of the preceding proposition, provides an estimate on certain Jacobi vector fields, with explicit dependence on $\|R\|$.

\begin{lemma} \label{lemma_Jacobi} Assume given the following data:
\begin{itemize}
\item $(X',g)$ asymptotically conical, a boundary-defining function $x$ and $V\in \nu_{sc}(X')$.
\item An estimate on the curvature tensor
$$\|R(x,\xi)\|\leq \rho^2(x)\leq c\,x^2$$
where $\rho(x)$ is strictly positive for $x>0$ and is independent of $\xi\in\partial X'$.

Notice that, by lemma \ref{lemma_local}, such an estimate always exists.
\item $\phi:X\longrightarrow X'$ a b-immersion, $y\in\partial X$ and $Z\in\nu_{sc}(X)$ such that $Z$ is non-zero along the curve $(x,y)\subseteq \Sigma$.
\end{itemize}

Let $J(s)$ denote the Jacobi vector field (depending on the parameter $x$) solution to
$$\nabla_s\nabla_s J(s)+R(J(s),\dot{\gamma})\dot{\gamma}=0,\ \ J(0)=Z[x,y],\ \nabla_sJ(0)=\nabla_Z V[x,y]$$
where $x\mapsto\gamma(x,s)$ denotes the curve of geodesics defined ($\forall s\in \R$) by $\gamma(x,0)=(x,y),\frac{\partial}{\partial\,s}\gamma(x,0)=V[x,y]$ and with the notation $\dot{\gamma}:=\frac{\partial}{\partial\,s}\gamma$. 

Then
$$J(s)=Z[x,y]+s\nabla_zV[x,y]+Q(s)\ \ (s\leq 1)$$
where $Z,\nabla_zV$ are extended along $\gamma(s)$ by parallel transport and $\|Q(s)\|=O(\rho(x)x)$.
\end{lemma}
\textbf{Proof}: Since $\|R(J(s),\dot{\gamma})\dot{\gamma}\|\leq \|R\|\cdot\|J(s)\|\cdot\|V\|^2\leq cx^2\|J(s)\|\cdot \|V\|^2$, intuitively we expect $J(s)\rightarrow A(s)$, where $A(s)$ is the solution of the equation
$$\nabla_s\nabla_s A(s)=0,\ \ A(0)=Z[x,y],\ \nabla_sA(0)=\nabla_ZV[x,y]$$
The proof can be divided in several steps and partially follows [J] (theorem 4.5.2) and [BK] (lemma 6.3.7). 

For simplicity, we will assume $\|V\|\leq 1$, $c\leq 1$.

\ 

\textbf{Step 1} Assume $g(s)$ is a $C^2$ function which, for some fixed $\eta\in \R$, satisfies
$$\ddot{g}+\eta\,g\leq 0,\ \ g(0)=\dot{g}(0)=0$$
Then
\begin{itemize}
\item if $\eta\leq 0$: $g(s)\leq 0$,\ \ $\forall s\geq 0$
\item if $\eta>0$: $g(s)\leq 0$,\ \ $0\leq s\leq \frac{\pi}{\sqrt{\eta}}$
\end{itemize}
Proof: Set

$f_\eta(s):=\left\{ \begin{array}{lr}
\frac{1}{\sqrt{\eta}}sin(\sqrt{\eta}\,s) & \eta>0\\
s & \eta=0\\
\frac{1}{\sqrt{-\eta}}sinh(\sqrt{-\eta}\,s) & \eta<0
\end{array} \right.$

\ 

ie, $f_\eta$ is the solution of the ODE
$$\ddot{f}+\eta\,f=0,\ \ f(0)=0,\ \dot{f}(0)=1$$
Since $g(s)$ has a second-order zero in $s=0$ while $f_\eta(s)$ has a first-order zero, we find $\frac{g}{f_\eta}_{|s=0}=0$.

Now set $d:=g'f_\eta-g\,f'_\eta$, so that $(\frac{g}{f_\eta})'=\frac{d}{f^2_\eta}$.

Notice that $d(0)=0$ and that 
$$d'(s)=g''f_\eta-g\,f_\eta''=(g''+\eta\,g)\,f_\eta\leq 0 \Leftrightarrow f_\eta\geq 0$$
Assume $\eta\leq 0$. Then $f_\eta\geq 0\Leftrightarrow s\geq 0$, so $d\leq 0$ for these $s$. This shows that $(\frac{g}{f_\eta})'\leq 0$, so $g\leq 0$.

The case $\eta>0$ is similar.

\ 

\textbf{Step 2} Consider the solution $b(s)$ of the ODE
$$\ddot{b}(s)=\|R(s)\|\cdot\|J(s)\|,\ \ b(0)=\dot{b}(0)=0$$
We want to prove that $\|J-A\|(s)\leq b(s)$. Clearly, it is enough to prove that, for all $P$ unit parallel vector fields along $\gamma$,
$$(J-A,P)-b\leq 0$$
Notice that, setting $\beta:=(J-A,P)$, $\beta$ is $C^2$ and satisfies
$$\ddot{\beta}\leq \|R(s)\|\cdot\|J(s)\|,\ \ \beta(0)=\dot{\beta}(0)=0$$
so that $g(s):=\beta-b$ satisfies the hypothesis of step 1 with $\eta=0$. Thus $\beta(s)\leq b(s)$, for $s\geq 0$.

\ 

\textbf{Step 3} From step 2, we also get
$$\ddot{b}\leq \|R(s)\|\cdot \|J-A\|+\|R(s)\|\cdot \|A(s)\|\leq \|R(s)\|\,b(s)+\|R(s)\|\cdot \|A(s)\|$$
Using the fact that $\|V\|$ is bounded and the uniformity of our estimate on $\|R\|$, we find $$\|R(s)\|\leq c\,\rho^2(x)\ \ (s\leq 1)$$
For simplicity, we will assume $c\leq 1$.

Notice that, since $\ddot{b}(s)\geq 0$, $b$ is non-negative for $s>0$. Thus
$$\ddot{b}\leq \rho^2(x)b(s)+\rho^2(x)\|A(s)\|\ \ (0\leq s\leq 1)$$
Let $a(s)$ be the solution to the ODE
$$\ddot{a}(s)=\rho^2(x)a(s)+\rho^2(x)\|A(s)\|\ \ a(0)=\dot{a}(0)=0$$
Then $g(s):=b-a$ satisfies the hypothesis of step 1 with $\eta=-\rho^2$, so $b(s)\leq a(s)$ and $\|J-A\|(s)\leq a(s)$, for $0\leq s\leq 1$.

\ 

\textbf{Step 4} Notice that $A(s)=Z[x,y]+s\nabla_ZV[x,y]$ where $Z[x,y]$ and $\nabla_ZV[x,y]$ are defined along $\gamma$ by parallel transport. 

To finish the proof, we now need to estimate $a(s)$.

Notice that $a$ is $C^2$ and is defined on some maximal interval $(\alpha,\beta)$.

Since $\ddot{a}(0)=\rho^2(x)\cdot\|Z\|>0$, $\dot{a}(s)$ is positive on some maximal connected interval $(0,\delta)$.

Assume $\delta<\beta$; in particular, this implies $\dot{a}(\delta)=0$. When $s\in (0,\delta)$, $a$ is monotone so $a(s)>0$ and $\ddot{a}(s)> 0$. Thus $\dot{a}(s)$ is non-decreasing, which contradicts $\dot{a}(\delta)=0$.

This shows that $\delta=\beta$, i.e. $\dot{a}(s)$ and $a(s)$ are positive for all $0<s<\beta$.

By hypothesis (and using $\|\nabla_ZV\|=O(x)$), $\ddot{a}\leq c(\rho^2a+\rho^2+s\rho^2x)$. Thus
$$2\ddot{a}\dot{a}\leq 2c\rho^2a\dot{a}+2c\rho^2\dot{a}+2c\rho^2xs\dot{a}$$
Integrating in $ds$, we find 
$$\dot{a}^2\leq c\rho^2a^2+2c\rho^2a+2c\rho^2xsa-2c\rho^2x\int a\,ds$$
By definition, $-c\rho^2a\leq -\ddot{a}+c\rho^2+c\rho^2xs$, so
$$-2c\rho^2x\int a\,ds\leq -2x\dot{a}+2c\rho^2xs+c\rho^2x^2s^2$$
Applying this above and using $c\leq c^2$, $\rho\leq x$, we find
\begin{eqnarray*}
\dot{a}^2 &\leq & c^2\rho^2a^2+2c\rho^2a+2c^2\rho^2xsa-2x\dot{a}+2c\rho^2xs+c^2\rho^2x^2s^2\\
(\dot{a}+x)^2 &\leq & (c\rho a+c\rho xs)^2+2c\rho xa+2c\rho x^2s+x^2\\
\dot{a}+x &\leq & c\rho a+c\rho xs+x
\end{eqnarray*}
When $s\leq 1$, this shows that $\dot{a}\leq c\rho a+c\rho x$ so, using Gronwall's inequality,
$$a(s)\leq c\rho x e^{c\rho s}\leq c_1\rho(x)x\ \ \ (s\leq 1)$$
In particular, $a(s)$ is defined up to $s=1$.
\qed
\section{Deformations of asymptotically conical SL submanifolds} \label{section_SLdefs}
Let $(X',g)$ be a fixed, 2n-dimensional, asymptotically conical manifold with a CY structure $(g,J,\omega,\Omega=\alpha+i\beta)$ on its interior $N:=X'\setminus \partial X'$. We will assume that $\omega\in \mbox{$^{sc}\Lambda^2$}X',\beta\in \mbox{$^{sc}\Lambda^n$}X'$; all these conditions are, for example, verified by $N=\C^n$ with its standard structures.

Let $\Sigma$ be a manifold with ends.

\begin{definition} An ``asymptotically conical special Lagrangian immersion'' (AC SL) of $\Sigma$ into $N$ is a (smooth) b-immersion $\phi:X\longrightarrow X'$ such that $\phi_{|\Sigma}$ is SL, where $X$ is some compactification of $\Sigma$.
\end{definition}

We refer to section \ref{section_SL} for examples of AC SL immersions.

The goal of this section is to study the integrability of ``infinitesimal SL deformations'' of a given AC SL $(\Sigma,\phi)$; as we will see further on, these are the normal vector fields on $\Sigma$ corresponding, via the isomorphism $(T\Sigma)^\perp\simeq\Lambda^1(\Sigma)$, to 1-forms $\nu$ on $\Sigma$ such that $(d\oplus d^*_g)\nu=0$.

Of course, it is important to specify ``boundary conditions'' for the allowed deformations. Using weighted Sobolev spaces, we will impose that our vector fields decay at infinity; i.e., we study deformations that ``keep infinity fixed''. Notice that the results of section \ref{section_immersions} show that such deformations automatically preserve the ``AC'' condition.

\begin{theorem} \label{theorem_ACSL1} Let $\phi: \Sigma \hookrightarrow N$ be an AC SL immersion. Then, for large $k$, small $\epsilon>0$ and using the isometry $\Lambda^1(\Sigma)\simeq T\Sigma^\perp$, the set
$$\mathcal{D}ef_{SL}(\Sigma,\phi):=\{V \in x^\epsilon H^k_b(\Lambda^1,sc): exp\,V \circ \phi \mbox{ is SL}\}$$
is smooth near the origin and has dimension dim $\mathcal{K}^1_\epsilon(\Sigma,\phi^*g)=\mbox{dim } H^1(\Sigma)+\mbox{dim } \mathcal{H}^0_{-1+\epsilon}(\Sigma)-1$: it thus depends both on the topology of $\Sigma$ and on the analytic properties of the link of $(\Sigma,\phi^*g)$.
\end{theorem}
\textbf{Proof}: To simplify the notation, we will identify $\Sigma$ with its image in $M$, omitting the dependence on the fixed immersion $\phi$. Consider the map
$$\begin{array}{rclcl}
F_c: \Lambda^1_c(\Sigma) & \simeq & \Lambda^0_c(T\Sigma^\perp) & \longrightarrow & \Lambda^2_c(\Sigma) \oplus C^{\infty}_c(\Sigma)\\
\nu & \simeq & V & \mapsto & (exp\,V)^*\omega_{|\Sigma} \oplus *_g[(exp\,V)^*\beta_{|\Sigma}]
\end{array}$$
where $exp$ is the exponential map of $(N,g)$ and $*_g$ is the Hodge star operator of $(\Sigma,\phi^*g)$.

Notice that, if $V\in x^\epsilon H^k_b(T\Sigma^\perp,sc)$ and using lemma \ref{lemma_Jacobi},
\begin{eqnarray*}
(exp\,V^*\omega)[p](e_i,e_j) &=& \omega[expV](expV_*e_i,expV_*e_j)\\
&=& \omega[expV](e_i+\nabla_{e_i}V+Q_i(1),e_j+\nabla_{e_j}V+Q_j(1))
\end{eqnarray*}
where $\{e_i\}$ is a local orthonormal basis for $(\Sigma,\phi^*g)$.

Using $\nabla\omega=0$, $\omega\in\mbox{$^{sc}\Lambda^2$}X'$ and the fact that $\Sigma$ is Lagrangian, the above terms may be examined one by one. For example:
\begin{itemize}
\item $\omega[expV](e_i,e_j)=\omega[p](e_i,e_j)=0$
\item $\omega[expV](e_i,\nabla_{e_j}^TV)=\omega[p](e_i,\nabla^T_{e_j}V)=0$
\item $|\omega[expV](e_i,\nabla_{e_j}^\perp V)|=|\omega[p](e_i,\nabla_{e_j}^\perp V)|\leq |\nabla_{e^j}^\perp V|$
\item $|\omega[expV](e_i,Q_j(1))|=O(x^2)$
\end{itemize}
This proves that $exp\,V^*\omega_{|\Sigma}\in x^{\epsilon+1}L^2_b(\Lambda^2\Sigma,sc)$.

More generally, one can prove that $F_c$ has an extension
$$F: x^\epsilon H^k_b(\Lambda^1,sc) \longrightarrow x^{\epsilon+1} H^{k-1}_b(\Lambda^2,sc) \oplus x^{\epsilon+1} H^{k-1}_b(\Sigma)$$
We want to prove that $F^{-1}(0)$ is smooth using the implicit function theorem. To be able to apply the Banach-space version of this theorem, it is sufficient to prove that $dF[0]$ is surjective and has finite-dimensional kernel.

Notice that
$$\frac{d}{dt}_{|t=0}exp(tV)^*\omega_{|\Sigma}=\mathcal{L}_V \omega_{|\Sigma}=(i_V\,d+d\,i_V)\omega_{|\Sigma} =d\,\nu$$
Linear algebra shows that $*_g i_V\beta_{|\Sigma}=(-1)^{n-1}\nu$; thus,
$$\frac{d}{dt}_{|t=0}*_g[exp(tV)^*\beta_{|\Sigma}]=*_g[\mathcal{L}_V\beta_{|\Sigma}]=(-1)^{n-1} d^*_g(*_gi_V\beta_{|\Sigma})=d^*_g\nu$$
so that $dF[0]=D_g$. It is thus Fredholm, as seen above. Surjectivity onto the whole space $x^{\epsilon+1} H^{k-1}_b(\Lambda^2,sc) \oplus x^{\epsilon+1} H^{k-1}_b(\Sigma)$ is actually false, but the following steps show that we may restrict our map to a smaller range and obtain surjectivity.
\begin{step} $Im(F_c) \subseteq Im (d_c) \oplus Im (d_c^*)$, i.e.:
$$\exists \sigma \in \Lambda^1_c(\Sigma): (exp\,V)^*\omega_{|\Sigma} = d_c \sigma,\ \ \ 
\exists \tau \in \Lambda^1_c(\Sigma): *_g[(exp\,V)^*\beta_{|\Sigma}] = d^*_c\tau$$
where $d_c$, $d_c^*$ denote the restrictions of $d$, $d^*$ to $\Lambda^*_c$.
\end{step} 
Proof: Extend $V$ to a compactly-supported vector field on $N$. Then $t\mapsto \phi_t:=exp(tV)$ is a $1$-parameter group of diffeomorphisms and 
\begin{eqnarray*}
\phi_t^*\omega & = & \phi_0^*\omega + \int_0^t \frac{d}{ds}(\phi_s^*\omega)\,ds\\
&=& \phi_0^*\omega +\int_0^t \phi^*_s(\mathcal{L}_{V}\,\omega)\,ds\\
&=& \phi_0^*\omega +\int_0^t d \, (\phi^*_s i_{V}\,\omega) \, ds\\
&=& \phi_0^*\omega +d(\int_0^t \phi^*_s i_{V}\,\omega\,ds)
\end{eqnarray*}
If we set $\sigma:=\int_0^1 \phi^*_s i_V\omega\,ds$, it is thus clear that $\sigma_{|\Sigma} \in \Lambda^1_c(\Sigma)$ and that $\phi^*\omega_{|\Sigma}=d\,(\sigma_{|\Sigma})$, as claimed. 

The second claim may be restated as
$$\exists \tau \in \Lambda^{n-1}_c(\Sigma): (expV)^*\beta_{|\Sigma} = d_c\tau$$
and can thus be proved as above.
\begin{step} $Im(d) \oplus Im(d^*) = Im (D_g)$, i.e.:
$$\forall \alpha, \beta \in x^\epsilon H^k_b(\Lambda^1,sc) \,\  \exists \gamma \in x^\epsilon H^k_b(\Lambda^1,sc): d\,\alpha \oplus d^*\beta = d \,\gamma \oplus d^* \gamma$$
\end{step}
Proof: Let $\gamma:=\alpha + d\,f$, for some $f \in x^{\epsilon-1}H^{k+1}_b(\Sigma)$. It is enough to show that $f$ can be determined so that $d^*\beta=d^*\alpha+d^*df$, i.e. $d^*(\beta-\alpha)=\Delta_g f$. So, it is enough to show that
$$Im(d^*:x^\epsilon H^k_b(\Lambda^1,sc) \longrightarrow x^{\epsilon+1} H^{k-1}_b(\Sigma))=Im(\Delta_g:x^{\epsilon-1}H^{k+1}_b \longrightarrow x^{\epsilon+1}H^{k-1}_b)$$

Notice that $Im(d^*) \supseteq Im (\Delta_g)$. Since $\epsilon <n-1$, $\Delta_g$ is surjective (on functions) so the above is necessarily true.
\begin{step} $Im (F) \subseteq Im (D_g)$
\end{step}
Proof:
$$Im(F) \subseteq \overline{Im(F_c)} \subseteq \overline{Im(d_c) \oplus Im (d^*_c)} \subseteq \overline{Im(d) \oplus Im (d^*)}=\overline{Im(D_g)}=Im(D_g)$$

\ 

The above proves that the restriction
$$F:x^\epsilon H^k_b(\Lambda^1,sc) \longrightarrow Im(D_g)$$
is a well-defined map between Banach spaces, such that $dF[0]$ has finite-dimensional kernel and is surjective. We may thus apply the implicit function theorem, proving that $F^{-1}(0)$ is smooth near the origin and has dimension dim $\mathcal{K}^1_\epsilon(\Sigma,\phi^*g)$, which we calculated in section \ref{subsection_1forms}. 
\qed
Recall that each SL immersion $exp\,V\circ\phi$ is, in particular, minimal; thus, standard regularity results prove that $exp\,V\circ\phi$ is smooth on $\Sigma$; i.e., $exp\,V\circ\phi\in C^\infty(\Sigma,N)$. Since $exp$ is a diffeomorphism on normal vector fields, this also shows that $V$ is smooth on $\Sigma$. 

\ 

\textbf{Remarks:}
\begin{enumerate}
\item Given $p\in \Sigma$, recall from section \ref{section_SL} the isometries
$$(\alpha\in T_p^*\Sigma) \stackrel{g}{\simeq} (W\in T_p\Sigma) \stackrel{J}{\simeq} (V \in T_p\Sigma^\perp)$$
The bundles $T^*\Sigma$, $T\Sigma$ over $\Sigma$ have natural extensions to bundles $^{sc}T^*X$, $^{sc}TX$ over $X$. It is using this fact that we can define operators $\mbox{Diff}^k_b$ and spaces $x^\delta H^k_b(\Lambda^1,sc)$, $x^\delta H^k_b(T\Sigma,sc)$.

Since $g$ is not defined on $\partial X$, however, it is not clear, a priori, if $T\Sigma^\perp$ has an extension up to $\partial X$. On the other hand, the Lagrangian condition shows that $T\Sigma^\perp\simeq T\Sigma$; thus, up to this identification, $^{sc}TX$ provides an extension of $T\Sigma^\perp$, allowing us to define the spaces $x^\delta H^k_b(T\Sigma^\perp,sc)$, used in the above proof.
\item Let $\alpha\in x^\epsilon H^k_b(\Lambda^1,sc)$ for $k$ large and $\epsilon>0$. By the Sobolev embedding theorem, notice that $\alpha$ is $C^1$ on $\Sigma$ and that $\|\alpha\|=O(x^\epsilon)$. Examining the formula for the induced metric on 1-forms shows, however, that this is not enough to ensure that $\alpha$ has a continuous extension up to $\partial X$; this corresponds to the fact that, in general, $\alpha\in \nu_{sc}^*(X)$ does not extend up to $\partial X$.

On the other hand, let $W \in x^\epsilon H^k_b(T\Sigma,sc)$ be the corresponding tangent vector field. Then:
\begin{itemize}
\item As above, $W$ is $C^1$ on $\Sigma$.
\item If we write $W=a\,\partial_x+b_i\,\partial_{y^i}=\frac{a}{x^2}(x^2\partial_x)+\frac{b_i}{x}(x\,\partial_{y^i})$, the fact $\|W\|_{ac}=O(x^\epsilon)$ implies that $\frac{a}{x^2}\rightarrow 0$, $\frac{b_i}{x}\rightarrow 0$.

Thus $W$ admits a continuous extension to zero on $\partial X$.
\item Notice that 
\begin{eqnarray*}
\nabla^T_{x\partial_x}W&=&x(\partial_x a)\partial_x+xa\,\Gamma_{xx}^x\partial_x+xa\,\Gamma_{xx}^j\partial_{y^j}+x(\partial_x{b_j})\partial_{y^j}\\
&&+xb_i\,\Gamma_{xi}^x\partial_x+xb_i\,\Gamma_{xi}^j\partial_{y^j}
\end{eqnarray*}
The fact that $\|\nabla_{x\partial_x}^TW\|_{ac}=O(x^\epsilon)$ thus implies that $\partial_xa,\partial_xb_j$ have extensions to zero on $\partial X$. 

In a similar way we can prove that $W$ has a $C^1$ extension to zero on $\partial X$.
\end{itemize}
The same holds for the normal field $V:=JW$. If we let $\nu_{1,b}$, $\nu_{1,sc}$ denote the $C^1$ analogues of $\nu_b$, $\nu_{sc}$, we thus find that $W\in\nu_{1,sc}(X)$, $V\in \nu_{1,sc}(X,X')$.
\item Although from a geometric point of view one is interested only in the smoothness over $\Sigma$ of the SL deformations $exp\,V\circ\phi$, it is important to point out that, under compactification, these deformed submanifolds will not, in general, be smooth up to $\partial X$. This is a standard situation in b-geometry: b-ellipticity is not sufficient to ensure regularity up to $\partial X$. 

One can however prove the existence of a certain asymptotic expansion of $V$ at $\partial X$, i.e. ``polyhomogeneity'' in the sense of [M1]. In our situation, this implies that $exp\,V\circ\phi\in C^1(X;X')\bigcap C^\infty(\Sigma;N)$. In any case, the results of section \ref{section_immersions} apply to show that $exp\,V\circ\phi$ is ($C^1$) asymptotically conical.
\item In this article we have given a purely metric characterization of the condition ``asymptotically conical''; this has the advantage of allowing for fairly general ambient spaces. When $(N,g)=(\C^n,g_{std})$ and the submanifold $\Sigma$ is minimal (as in our case), our definition implies alternative, ``set-approximating'' definitions, as follows. 

Recall that the components of a minimal immersion are harmonic functions of $(\Sigma,g)$. The results of section \ref{subsection_functions} thus apply to show that there exists a minimal cone, lying in $\C^n$, to which $\Sigma$ ``converges'', with speed $o(r)$. When $\Sigma$ is Lagrangian, the cone is SL. 
\item As mentioned in the introduction, the results of this paper are very close to those conjectured in [J1], [J2] and proved in [M]. The main difference between this work and that of Joyce and Marshall is in the set-up of the problem: while we measure our deformations with respect to the given initial AC SL submanifold (this ``intrinsic approach'' works in any AC CY), Joyce prefers to postulate the existence of a SL cone, then measures the deformations of AC SL submanifolds with respect to that cone. This ``extrinsic approach'' works only in $\C^n$ (where cones can be defined); in this ambient space, their results could be reconstructed from ours, starting with the observation in remark 4 (above) which (adopting the terminology of [J1] and [J2]) basically states that, in $\C^n$, any AC SL submanifold is ``weakly asymptotic'' to some cone.
\end{enumerate}

As a final step, it is interesting to understand the dependence of the above construction on the weight $\epsilon$: in theory, changing the weight changes the class of allowed deformations and, thus, the dimension of the space of AC SL deformations.

In particular, it is interesting to consider what happens with respect to the spaces $x^{\frac{n}{2}} H^k_b$, corresponding to $L^2$-decay of the deformations. As shown in the course of the proof of theorem \ref{theorem_ACSL1}, to get a smooth structure dependent on this weight one would need an extension of $F_c$ to a map
$$F:x^{\frac{n}{2}} H^k_b(\Lambda^1,sc) \longrightarrow x^{\frac{n}{2}+1} H^{k-1}_b(\Lambda^2,sc)\oplus  x^{\frac{n}{2}+1} H^{k-1}_b(\Sigma)$$
However, the existence of such an extension depends on the properties of $exp$. According to lemma \ref{lemma_Jacobi}, the curvature $R$ of $N$ introduces a perturbation $Q$ which depends on the rate of decay of $\|R\|$; if $\|R\|$ does not decay sufficiently fast (with respect to the weights in consideration), the above extension does not exist, because the perturbation is too big with respect to the weights.

Consider, for example, the case $(N,g)=(\C^n,g_{std})$. In this case, in lemma \ref{lemma_Jacobi}, $Q\equiv 0$ and the extension exists. The whole proof of the theorem carries through as before, yielding a second, smaller, set of deformations.

We can thus prove the following theorem, which actually holds for any $(N,g)$ with sufficiently fast curvature decay.
\begin{theorem} \label{theorem_ACSL2} If $(N,g)=(\C^n,g_{std})$, the set
$$\mathcal{D}ef_{SL}^{L^2}(\Sigma,\phi):=\{V\in x^\frac{n}{2} H^k_b(\Lambda^1,sc): exp\,V \circ \phi \mbox{ is SL}\}$$
is a smooth submanifold of $\mathcal{D}ef_{SL}(\Sigma,\phi)$.

It has dimension dim $\mathcal{K}^1_\frac{n}{2}(\Sigma,\phi^*g)=\mbox{dim } H^1_c(\Sigma)$ and thus depends only on the topology of $\Sigma$.
\end{theorem}

\

\mbox{\Large\textbf{Bibliography}}
\begin{description}
\item[\mbox{[A]}] Aubin, T., Some Nonlinear Problems in Riemannian Geometry, Springer, Berlin, 1998
\item[\mbox{[B]}] Bartnik, R., The mass of an Asymptotically Flat Manifold, Comm. Pure and Applied Math., vol. XXXIX (1986), pp. 661-693
\item[\mbox{[BK]}] Buser, P. and Karcher, H., Gromov's almost flat manifolds, Asterisque 81, 1981
\item[\mbox{[C]}] Calabi, E., M\'{e}triques Kaehl\'{e}riennes et fibr\'{e}s holomorphes, Ann. scient. \'{E}c. Norm. Sup. (4), 12 (1979), pp. 269-294
\item[\mbox{[CZ]}] Christiansen, T. and Zworski, M., Harmonic functions of polynomial growth on certain complete manifolds, Geom. and Func. Analysis, vol. 6, n. 4 (1996), pp. 619-627
\item[\mbox{[G]}] Goldstein, E., Minimal Lagrangian tori in Kaehler Einstein manifolds, math.DG/0007135 (pre-print)
\item[\mbox{[HL]}] Harvey, R. and Lawson, H.B., Calibrated geometries, Acta Mathematica, 148 (1982), pp. 47-156
\item[\mbox{[H]}] Haskins, M., Special Lagrangian Cones, math.DG/0005164 (pre-print)
\item[\mbox{[J]}] Jost, J., Riemannian geometry and geometric analysis, Springer, 1998
\item[\mbox{[J1]}] Joyce, D., On counting special Lagrangian homology 3-spheres, hep-th/9907013 (pre-print)
\item[\mbox{[J2]}] Joyce, D., Lectures on Calabi-Yau and special Lagrangian geometry, part I of Gross, M. (et al.), Calabi-Yau Manifolds and Related Geometries, Springer-Verlag, 2003.
\item[\mbox{[M]}] Marshall, S., Deformations of special Lagrangian submanifolds, Oxford D.Phil thesis, 2002.
\item[\mbox{[ML]}] McLean, R., Deformations of calibrated submanifolds, Comm. An. and Geom., 6 (1998), pp. 705-747
\item[\mbox{[MO]}] McOwen, R., The Bahavior of the Laplacian on Weighted Sobolev Spaces, Comm. Pure and Applied Math., vol. XXXII (1979), pp. 783-795
\item[\mbox{[M1]}] Melrose, R., The Atiyah-Patodi-Singer Index Theorem, A.K.Peters, Wellesley, MA, 1994
\item[\mbox{[M2]}] Melrose, R., Spectral and Scattering Theory for the Laplacian on Asymptotically Euclidean Spaces, in Spectral and Scattering Theory, Ikawa M. ed., lecture notes in pure and applied mathematics, vol. 161 (1994), Marcel Dekker Inc.
\item[\mbox{[NW]}] Nirenberg, L., and Walker, H., The null spaces of elliptic partial differential operators in $\R^n$, J. Math. Anal. Appl., 42 (1973), pp. 271-301
\item[\mbox{[Oh]}] Oh, Y.-G., Second variation and stabilities of minimal lagrangian submanifolds in Kaehler manifolds, Invent. math., 101 (1990), pp. 501-519
\item[\mbox{[P]}] Pacini, T., Flows and Deformations of Lagrangian Submanifolds in Kaehler-Einstein Geometry, Ph.D. Thesis, Univ. of Pisa (2002)
\item[\mbox{[Y]}] Yau, S.T., On the Ricci curvature of a compact Kaehler manifold and the complex Monge-Amp\`{e}re equation. I, Comm. Appl. Math., 31 (1978) n.3, pp. 339-411 
\end{description} 

\ 

\ 

Tommaso Pacini (Imperial College, London/University of Pisa)

\ 

Email: pacini@paley.dm.unipi.it

\ 

Subject class: 53C38; 58J05
\end{document}